\documentclass[twoside,11pt]{article}
\renewcommand{\baselinestretch}{1.5}

\usepackage{graphicx}
\usepackage{graphics,ifthen}
\usepackage{latexsym}
\usepackage{amsfonts}
\usepackage{amsmath}
\usepackage{amssymb}%
\usepackage{amsthm}
\usepackage{amscd}
\usepackage{natbib}
\begin{document}
\input{latexP.sty}
\input{referencesP.sty}
\input{epsf.sty}
\input{psfig.sty}

\def\Definition{\stepcounter{definitionN}
    \Demo{Definition\hskip\smallindent\thedefinitionN}}
\def\EndDefinition{\EndDemo}
\def\EndExample{\qed\EndDemo}
\def\Category#1{\centerline{\Heading #1}\rm}

\def\e{\text{\hskip1.5pt e}}
\newcommand{\eps}{\epsilon}
\newcommand{\remarks}{\noindent {\bf Remarks:\ }}
\newcommand{\note}{\noindent {\bf Note:\ }}
\newcommand{\Lower}[2]{\smash{\lower #1 \hbox{#2}}}
\newcommand{\ben}{\begin{enumerate}}
\newcommand{\een}{\end{enumerate}}
\newcommand{\bi}{\begin{itemize}}
\newcommand{\ei}{\end{itemize}}
\newcommand{\hp}{\hspace{.2in}}
\newtheorem{thm}{Theorem}[section]
\newtheorem{examp}{Example}[section]
\newtheorem{defin}{Definition}[section]
\newtheorem{prop}{Proposition}[section]
\newtheorem{lem}[thm]{Lemma}
\newtheorem{cor}[thm]{Corollary}
\newtheorem{rem}[thm]{Remark}
\newtheorem{algorithm}[thm]{Algorithm}
\newcommand{\rb}[1]{\raisebox{1.5ex}[0pt]{#1}}
\newcommand{\mc}{\multicolumn}

\def\V{\mathbb{V}}
\def\P{\mathbb{P}}
\def\PY{{\cal PD}}
\def\Z{{\bf Z}}
\def\D{{\bf D}}
\def\X{{\bf X}}
\def\U{{\bf U}}
\def\Y{{\bf Y}}
\def\B{\mbox{\boldmath{$\beta$}}}
\def\sB{\mbox{{\scriptsize\boldmath{$\beta$}}}}
\newcommand{\fC}{\cal C_{\S}}
\newcommand{\bp}{\mathbf p}
\newcommand{\bX}{\mathbf X}
\def\Var{\mathop{\rm Var}\nolimits}%
\def\Cov{\mathop{\rm Cov}\nolimits}%
\newcommand{\bone}[1]{{\Bbb I}_{\left\{#1\right\}}}

\def\R{\mathcal{R}}
\def\X{\mathbf{X}}
\def\XX{\mathcal{W}}
\def\VV{\mathbf{V}}
\def\UU{\mathbf{U}}
\def\s{\mathbf{s}}
\def\S{\mathbf{S}}
\def\C{\mathbb{C}}
\def\E{\mathbb{E}}
\def\c{\mathbf{c}}
\def\p{\mathbf{p}}
\def\np{n(\p)}
\def\Np{N(\p)}
\def\x{\mathbf{X}}
\def\betax{\mathbf{X}^{\prime}\mbox{\boldmath{$\beta$}}}
\def\sbetax{\mathbf{X}^{\prime}\mbox{\boldmath{$\beta$}}}
\def\tbetax{\mathbf{X}^{\prime}\mbox{\scriptsize\boldmath{$\beta$}}}
\def\bbeta{\mbox{\boldmath{$\beta$}}}
\def\sbbeta{\mbox{\boldmath{$\beta$}}}
\def\tbbeta{\mbox{\scriptsize\boldmath{$\beta$}}}
\def\T{\mathbf{T}}
\def\D{\mathbf{D}}
\def\Q{\mathbf{Q}}
\def\XI{\mbox{\boldmath{$\Omega$}}}
\def\y{\mathbf{y}}
\def\z{\mathbf{z}}
\def\W{\mathbf{W}}
\def\Y{\mathbf{Y}}
\def\Z{\mathbf{Z}}
\def\vv{\mathbf{v}}
\newcommand{\I}[1]{\,\mathbb{I}(#1)\,}
\def\CP{\mathcal{P}}
\def\CM{\mathcal{M}}

\def\Report{Unimodal Density}
\def\Author{Ho}
\pagestyle{myheadings} \markboth{\Author}{\Report}
\thispagestyle{empty} \bct\Heading Bayesian Nonparametric Estimation of a Unimodal Density\\
via two $\S$-paths\footnote{
\indent\eightit Key words and phrases.
                \rm
Species sampling model, Species sampling mixture model, 
Rao--Blackwellization, Sequential importance sampling, Markov chain
Monte Carlo,
Accelerated path sampler.
}\lbk\lbk \smc Man-Wai Ho
\rm\lbk
(\today) \ect 
\begin{abstract}
A Bayesian nonparametric method for unimodal densities on the real
line is provided by considering a class of species sampling mixture
models containing random densities that are unimodal and not
necessarily symmetric. This class of densities generalize the model
considered by Brunner~(1992), in which the Dirichlet process is
replaced by a more general class of species sampling models. A novel
and explicit characterization of the posterior distribution via a
finite mixture of two dependent $\S$-paths is derived. This results
in a closed-form and tractable Bayes estimator for any unimodal
density in terms of a finite sum over two $\S$-paths. To approximate
this class of estimates, we propose a sequential importance sampling
algorithm that exploits the idea of the accelerated path sampler, an
efficient path-sampling Markov chain Monte Carlo method. Numerical
simulations are given to demonstrate the practicality and the
effectiveness of our methodology.
\end{abstract}
\rm

A Bayesian nonparametric method for unimodal densities on the real
line is provided by considering a class of species sampling mixture
models containing random densities that are unimodal and not
necessarily symmetric. This class of densities generalize the model
considered by Brunner~(1992), in which the Dirichlet process is
replaced by a more general class of species sampling models. A novel
and explicit characterization of the posterior distribution via a
finite mixture of two dependent $\mathbf{S}$-paths is derived. This
results in a closed-form and tractable Bayes estimator for any
unimodal density in terms of a finite sum over two
$\mathbf{S}$-paths. To approximate this class of estimates, we
propose a sequential importance sampling algorithm that exploits the
idea of the accelerated path sampler, an efficient path-sampling
Markov chain Monte Carlo method. Numerical simulations are given to
demonstrate the practicality and the effectiveness of our
methodology.

\section{Introduction}\label{sec:intro}
Statisical theory usually assumes that data come from a distribution
that is symmetric and unimodal at zero, such as a normal
distribution or a Student's $t$ distribution. However, it is common
in real-life applications that underlying distribution of response
variable, even though unimodal, may not be symmetric about its mode
which is different from zero. For more information, see
Dharmadhikari and Joag-Dev (1988) and Bertin, Cuculescu and
Theodorescu (1997). There is a vast amount of literature on
nonparametric estimations of unimodal densities and the mode from a
frequentist viewpoint including early works of Granander~(1956),
Parzen~(1962), Chernoff~(1964), Robertson~(1967), Venter~(1967),
Prakasa Rao~(1969), Wegman~(1969, 1970a, 1970b, 1971), other further
studies by Lye and Martin~(1993), Bickel and Fan~(1996), Wang~(1996)
and Birg\'e~(1997) and among others. Some recent methods are, for
example, a recursive method in Cheng, Gasser, and Hall~(1999),
kernel-based methods in Hall and Huang~(2001, 2002), and other
parametric models in Fern\'{a}ndez and Steel~(1998), Jones~(2004)
and Ferreira and Steel~(2006).

From a Bayesian viewpoint, Brunner~(1992) gave a nonparametric
solution to the problems by assuming a mixture representation same
as that in~(\ref{mixture1}) wherein the mixing distribution $G$ is a
Dirichlet process~(Ferguson (1973)) for a unimodal density with a
general mode $\theta$ on the real line $\R$. The posterior
distribution and the Bayes estimate of the unimodal density can be
characterized in terms of random partitions~(see, e.g, Lo~(1984) and
Lo, Brunner and Chan~(1996) for these well-established theoretical
results on Dirichlet process mixture models).

This paper is devoted to improving Brunner's results and developing
an efficient numerical method for practical usage of the Bayes
solutions. A class of unimodal densities with mode at $\theta$ of
interest is defined by
\begin{equation}\label{mixture1}
    f(t|G,\theta) = \int
\frac{1}{X} \left[\I{0 < t-\theta \leq X}-\I{X \leq t-\theta <
0}\right] G(dX),\qquad t \in \R,
\end{equation}
where $\I{B}$ is the indicator of an event $B$ and $G$ is from the
class of species sampling models developed in Pitman~(1995, 1996),
of which the Dirichlet process is a member. All the results follow
are therefore applicable to Brunner's model as his model is a
special case of~(\ref{mixture1}). The validity of the mixture
representation for all unimodal densities given in the right hand
of~(\ref{mixture1}) can be justified by noting equality between its
integral when $\theta=0$ and the distribution function of any
unimodal density with mode at zero given in Feller~(1971, page~158).

The posterior distribution of~(\ref{mixture1}), like Brunner's
model, can also be characterized in terms of random partitions, as
the models are special cases of the species sampling mixture model
defined in Ishwaran and James~(2003) which takes the same form
as~(\ref{mixture1}) with the kernel ${X}^{-1} \[\I{0 < t-\theta \leq
X}-\I{X \leq t-\theta < 0}\]$ replaced by any density function in
$t$ given $\theta$ and $X$. In this work, by utilizing the special
and nice features of the kernel
in~(\ref{mixture1})~(see~(\ref{twocases})) and noticing irrelevancy
of some information carried by a partition in characterizing the
posterior distribution, we are able to refine the partition-based
results to show that the unimodal densities possess special
structures related to two $\S$-paths, where an $\S$-path is a random
vector defined in Brunner and Lo~(1989)~(see also Dykstra and
Laud~(1981). Generally speaking, there exists a tractable
characterization of the posterior distribution via some
combinatorial structures that are considerably less complex than
partitions. Such a characterization is known to be the first
explicit type that is based on two $\S$-paths. Similar phenomena
based on one single $\S$-path could be found in Bayes estimations of
symmetric unimodal or decreasing densities by Brunner and Lo~(1989),
Brunner~(1995) and Ho~(2006b) and monotone hazard functions by
Dykstra and Laud~(1981), Lo and Weng~(1989) and Ho~(2006a), as the
space of $\S$-paths is considerably smaller than that of
partitions~(Brunner and Lo 1989). Intuitively, this characterization
depending on two $\S$-paths can be explained by the fact that there
are two~(possibly different) non-increasing curves on each side of
the mode in unimodal densities, but not only one~(identical on
either sides) in symmetric unimodal densities of which can be
characterized in terms of one $\S$-path~(Albert~Y.~Lo, private
conversation).

It is recognized that if one could efficiently sample the two
$\S$-paths in this context, this would lead to more parsimonious
methods for inference. Motivated by the co-existence of and the
resemblance in constructions of an SIS algorithm and a Gibbs sampler
for sampling random partitions in many Bayesian mixture models~(Lo,
Brunner and Chan~(1996) and Ishwaran and James~(2003)), we
propose~(in Section~\ref{sec:SIS}) a novel sequential importance
sampling~(SIS) method~(Kong, Liu and Wong~(1994) and Liu and
Chen~(1998)), dubbed \textit{sequential importance path~(SIP)
sampler}, for sampling directly one single $\S$-path in the
aforementioned models under monotonicity constraints by borrowing
the idea behind the success of an efficient Markov chain Monte
Carlo~(MCMC) method introduced in Ho~(2002, 2006a, 2006b) that
serves the same purpose. Then, a natural SIS scheme based on
applications of the SIP sampler is introduced for sampling the
unknown mode $\theta$ and the two $\S$-paths in
evaluating/approximating posterior quantities for models
in~(\ref{mixture1}).

\subsection{Some backgrounds on species sampling
models}\label{ssmodel} Pitman~(1995, 1996) developed the class of
species sampling models that corresponds to the set of all random
probability measure of the form
\begin{equation}\label{randomP}
P(\cdot) = \sum_{k} W_k \delta_{V_k}(\cdot) +
\left(1-\sum_{k}W_k\right) H(\cdot),
\end{equation}
where $0<W_k<1$ are random weights such that $\sum_{k}W_k \leq 1$,
independently of $V_k$, which are i.i.d. random variables with some
non-atomic distribution $H$ over a measurable Polish space, and
$\delta_{V_k}(\cdot)$ is a Dirac probability measure degenerate at
$V_k$. This class includes a large number of well-known random
processes, for instance, the Dirichlet process and its two-parameter
extension, called the two-parameter Poisson-Dirichlet
process~(Pitman and Yor~(1997)), the class of finite-dimensional
Dirichlet priors discussed in detail in Ishwaran and
Zarepour~(2002a, 2002b), and the homogeneous normalized random
measures with independent increments discussed in Regazzini, Lijoi
and Pr\"unster~(2003).

Suppose $\X = (X_1,\ldots,X_{N})$ is a random sample
from~\mref{randomP}. The joint marginal distribution of $\X$ is
determined by the prediction rule, $\Pr\left\{X_1 \in \cdot \right\}
= H(\cdot)$ and
\begin{equation}
\label{gBM}
    \Pr\left\{X_{k+1} \in \cdot | X_1, \ldots, X_{k}\right\}
    = \ell_{0,k} \,H(\cdot) + \sum_{j=1}^{N_k}
    \ell_{j,k} \, \delta_{X^{\ast}_j} (\cdot), \quad
    k=2,\ldots,N-1,
\end{equation}
where $H$ is non-atomic and $\ell_{0,k}$ and $\ell_{j,k}$ are
non-negative measurable functions of $X_1,\ldots,X_k$. The above
prediction rule conveys that given $X_1, \ldots, X_{k}$, which
correspond to $N_k$ unique values $X_1^{\ast}, \ldots,
X_{N_{k}}^{\ast}$ of respective numbers of duplicates $e_1, \ldots,
e_{N_{k}}$, then the next observation $X_{k+1}$ takes the same value
as $X_j^{\ast}$ with probability $\ell_{j,k}$, $j = 1,\ldots,N_k$;
otherwise it takes a new value from $H$ with probability
$\ell_{0,k}$. As a consequence of the
exchangeability of $(X_1,\ldots,X_N)$
, Pitman~(1996) shows that the distribution of $X_1, \ldots, X_N$,
denoted by $\mu(d\X)$, is uniquely characterized by the joint law of
its unique values and an exchangeable partition probability
function~(EPPF)
\begin{equation}\label{symmetric}
\pi(\p) = \chi(e_1,\ldots,e_{\Np})
\end{equation}
induced by the unique values. That is,
$$
\mu(d\X) =\pi(\p)\prod_{k=1}^{\Np}H(dX_k^{\ast}),
$$
where $\p =\{C_1, \ldots, C_{\Np}\}$ of the integers
$\{1,\ldots,N\}$ is a partition of $\Np$ cells induced by $C_k =
\{j: X_j = X_k^{\ast}\}$ and $\chi$ is a unique symmetric function
depending only upon $e_k$, the number of elements in or the size of
$C_k$, $k=1,\ldots,N(\p)$~(see Pitman~(1996) and Ishwaran and
James~(2003, Section~2) for more information).

\section{A posterior distribution via $\S$-paths}
\label{sec:post}
This section concerns Bayes estimation of a unimodal density on the
line $\R$ with a general mode $\theta$, defined by the species
sampling mixture model in~(\ref{mixture1}). Suppose we observe $N$
i.i.d. observations $\T=(T_1, \ldots, T_N)$ from~(\ref{mixture1})
and assume any prior $\pi(d\theta)$ for $\theta$. Given $\T$, denote
$\CP(dG|\theta,\T)$ and $\CP(d\theta|\T)$ as the posterior
distribution of $G$ given $\theta$ and the posterior distribution of
$\theta$, respectively. The posterior distribution of the pair
$(G,\theta)$ in~(\ref{mixture1}) can always be determined by the
double expectation formula,
\begin{equation}\label{double}
\E[h(G,\theta) |\T ] = \E\{\E[h(G,\theta) |\theta,\T ]|\T\} =
\int_\R \int_{\CM} h(G,\theta) \CP(dG|\theta,\T) \CP(d\theta|\T),
\end{equation}
where $h$ is any nonnegative or integrable function and $\CM$ is the
space of probability measures over $\R$. Let us first look at
$\CP(dG|\theta,\T)$ and then discuss $\CP(d\theta|\T)$ later on.

Suppose $\theta$ is given. We can always assume that
\begin{equation}\label{obs}
(T_1-\theta,\ldots, T_N-\theta) = \Z \cup \Y = (Z_{N-n},Z_{N-n-1},
\ldots,Z_{1}) \cup (Y_1,Y_2,\ldots,Y_{n}),
\end{equation}
where $Z_{N-n}<Z_{N-n-1} < \cdots < Z_{1}<0$ and $0<Y_1 < Y_2 <
\cdots < Y_{n}$. Denote the missing variables in~(\ref{mixture1}) by
$\X = (X_1,\ldots,X_N)$. 
It is worthy of note that once an observation is taken
from~(\ref{mixture1}), the kernel can be well-simplified according
to two mutually exclusive situations, that is, the likelihood of an
observation $T_k$ in $\T$ is given by
\begin{equation}\label{twocases}
f(T_k|G,\theta) = \left\{
\begin{array}{lll}
\int (-X^{-1}) \I{X \leq T_k-\theta} G(dX) && T_k-\theta < 0 \\
\int X^{-1} \I{T_k-\theta \leq X} G(dX) && T_k-\theta > 0.
\end{array}
\right.
\end{equation}
The distinctiveness of the kernel yields a similar
simplification~(see~(\ref{jointp}) and~(\ref{varphi})) in the
posterior distribution of $G$ given $\theta$ in terms of partitions
$\p$ of the integers $\{1,\ldots,N\}$ in~(\ref{postp}), readily
available from Theorems~1 and~2 in Ishwaran and James~(2003). This
implies that the $n$ resulting positive observations after
subtraction of $\theta$, $Y_1,\ldots,Y_n$, can only ``cluster'' with
one another but not any negative observation or vice versa. Hence,
it is eligible to ``split'' the partition $\p$ of the $N$
integers/observations into two non-overlapping partitions $\p^{+}$
and $\p^{-}$. Write $\p = \p^+ \cup \p^-$. Without loss of
generality, assume that $\p^{+} = \{C_1,\ldots,C_{N(\p^+)}\}$
denotes the partition of the $n$ positive observations and $\p^{-} =
\{C_{N(\p^+)+1},\ldots,C_{\Np}\}$ of the remaining $N-n$ negative
observations. Define
\begin{equation}\label{chistar}
\pi(\p^-|\p^+) := 
\frac{\chi(e_1,\ldots,e_{N(\p^+)},e_{N(\p^+)+1},\ldots,e_{\Np})}
{\chi(e_1,\ldots,e_{N(\p^+)})}
= \frac{\pi(\p)}{\pi(\p^+)},
\end{equation}
where $\pi(\cdot)$ is defined in~(\ref{symmetric}), such that
\begin{equation}\label{symmetric2}
\pi(\p) = \pi(\p^-|\p^+) \times \pi(\p^+).
\end{equation}
These, together with the facts that the second line
of~(\ref{twocases}) resembles, while the other line is symmetrical
to, the scaled mixture of uniform representation of a symmetric
unimodal density with mode at zero due to Khintchine~(1938) and
Shepp~(1962), yield a posterior distribution of $G$ given $\theta$,
which is expressible in terms of two dependent $\S$-paths, as a
consequence of applications of Theorem~2.1 and Corollary~2.2 in
Ho~(2006b).

Let us fix some notation before stating the main results. Define an
integer-valued vector $\S = (S_0, S_1, \ldots, S_{n-1}, S_n)$
(Dykstra and Laud~(1981) and Brunner and Lo~(1989)), referred to as
an $\S$-path (of $n+1$ coordinates), which satisfies
(i) 
$S_0 = 0$ and $S_n = n$;
(ii) 
$S_j \leq j$, $j=1,\ldots,n-1$; and
(iii) 
$S_j \leq S_{j+1}$, $j=1, \ldots, n-1$.
A path $\S$ is said to \textit{correspond to} one or many partitions
$\p$ of the integers $\{1,\ldots,n\}$, provided that (i)~labels of
the maximal elements of the $\Np$ cells in $\p$ coincide with
locations $j$ at which $S_j
> S_{j-1}$, and (ii)~size $e_k$ of the cell $C_k$ for all
$k=1,\ldots,\Np$ with a maximal element $j$, $j=1,\ldots,n$, is
identical to $S_j - S_{j-1}$. Let ${\mathbb{C}}_{\S}$ denote the
collection of partitions that correspond to a given $\S$. Then, the
total number of partitions in ${\mathbb{C}}_{\S}$ is given
by~(Brunner and Lo~(1989))
\begin{equation}\label{card}
|{\mathbb{C}}_{\S}| = \displaystyle\prod_{j=1:S_j>S_{j-1}}^n
\binom{j - 1-S_{j-1}}{j-S_{j} }.
\end{equation}
See Ho~(2002) for more discussion of the relation between $\p$ and
$\S$. Following from the symmetric definition of $\chi$
in~(\ref{symmetric}), we have
\begin{equation}\label{priorS}
\pi(\p) = \chi(e_1,\ldots,e_{\Np}) := \chi(\mathcal{M}_{1,n}(\S)),
\qquad \mbox{if } \p \in {\C}_{\S},
\end{equation}
where, for any integer $1\leq a<b\leq n$,
$$
\mathcal{M}_{a,b}(\S) =
\{S_j-S_{j-1}:S_j>S_{j-1},j=a,a+1,\ldots,b\}.
$$
Write $\sum_{{\S}}$ as summing over all paths $\S$ of the same
number of coordinates, and $\prod_{\{j^{\ast}|\S\}}$ and
$\sum_{\{j^{\ast}|\S\}}$ as $\prod_{j=1:S_j>S_{j-1}}^{n}$ and
$\sum_{j=1:S_j>S_{j-1}}^{n}$ conditioning on $\S$, respectively.

\begin{thm} \label{Prop1}
{\rm Suppose $\theta$ is given and ${\T}$ are $N$ i.i.d.
observations from~(\ref{mixture1}). That is,~(\ref{obs}) holds.
Then, the distribution of $\X$ given $\theta$ and $\T$ can be
summarized by a joint law of $(\VV,\UU),(\S^-, \S^+)|\theta,\T$
defined as follows.
\begin{itemize}
\item[(i)] Given $(\theta,\T)$, two paths $\S^+ = (0, S_1^+, \ldots, S_{n-1}^+, n)$
and $\S^- = (0, S_1^-, \ldots, S_{N-n-1}^-,\break N-n)$ of $n+1$ and
$N-n+1$ coordinates, respectively, have a~(discrete) joint
distribution $\pi( \S^-,\S^+|\theta,\T) \propto
\phi^+_\theta({\S^+},\T) \times \phi^-_\theta({\S^-},\S^+,\T) $,
where
\begin{equation} \label{ZSphi}
    \phi^+_\theta( {\S}^+,\T) =
    \left|{\mathbb{C}}_{\S^+}\right|\chi(\mathcal{M}_{1,n}(\S^+))
    \prod_{\{j^{\ast}|\S^+\}}
    \int_{Y_j}^{\infty}
    U_j^{-(S_j^+-S_{j-1}^+)} H(dU_j)
\end{equation}
and
\begin{equation} \label{ZSphi2}
    \phi^-_\theta( {\S}^-,\S^+,\T) =
    \left|{\mathbb{C}}_{\S^-}\right|
    \frac{\chi(\mathcal{M}_{1,n}(\S^+),\mathcal{M}_{1,N-n}(\S^-))}
    {\chi(\mathcal{M}_{1,n}(\S^+))}
    \prod_{\{j^{\ast}|\S^-\}}
    \int_{-\infty}^{Z_j}
    (-V_j)^{-(S_j^--S_{j-1}^-)} H(dV_j)
\end{equation}
with $|\mathbb{C}_{\S}|$ and $\chi(\cdot)$ defined in~(\ref{card})
and~(\ref{priorS}), respectively.

\item[(ii)] Given $(\S^-,\S^+)$ and $(\theta,\T)$, there exist $N(\S^+) =
\sum_{j=1}^{n}\I{S_{j}^+>S_{j-1}^+}$ positive and $N(\S^-) =
\sum_{j=1}^{N-n}\I{S_{j}^->S_{j-1}^-}$ negative unique values on
$\R$ among $\{X_1,\ldots,X_N\}$, denoted by $\UU = \{U_j:
S_j^+>S_{j-1}^+,j=1,\ldots,n\}$ and $\VV = \{V_j:
S_j^->S_{j-1}^-,j=1,\ldots,N-n\}$, respectively. They are
distributed, conditionally independent of one another, as
\begin{equation} \label{alphaSj}
    H_j^+(dU_j \vert\S^+,\Y) \propto \I{Y_j \leq U_j} U_j^{-(S_j^+-S_{j-1}^+)}
    H(dU_j),
\end{equation}
and
\begin{equation}
\label{alphaSj2}
    H_j^-(dV_j \vert\S^-,\Z) \propto \I{V_j \leq Z_j} (-V_j)^{-(S_j^--S_{j-1}^-)}
    H(dV_j),
\end{equation}
respectively.
\end{itemize}
}
\end{thm}

\begin{thm}\label{Prop2} {\rm
For any nonnegative or integrable function $g$, the law of $G$ given
$\theta$ and $\T$ is characterized by
\begin{eqnarray}
&&\hspace*{-0.5in}\displaystyle\int_{\CM}g(G)\CP(dG|\theta,\T) \nonumber\\
&&\hspace*{-0.5in}\quad = \displaystyle \sum_{\S^+} \sum_{\S^-}
\left[\rule{0in}{5ex}
  \int_{\R^{N(\S^+)+N(\S^-)}} \left\{\int_\CM g(G)
  \CP(dG|\VV,\UU,\S^-, \S^+,\theta,\T)\right\}\right.\nonumber\\
&& \hspace*{0.5in} \left.\rule{0in}{5ex}
  \prod_{\{j^{\ast}|\S^-\}} H_j^-(dV_j\vert\S^-,\Z)
  \prod_{\{j^{\ast}|\S^+\}} H_j^+(dU_j\vert\S^+,\Y)\right]
  \pi(\S^-,\S^+|\theta,\T),\label{postS}
\end{eqnarray}
where $\CP(dG|\VV,\UU, \S^-,\S^+,\theta,\T)$ is equivalent in
distribution to $\CP(dG|\X^{\ast},\p,\theta,\T)$ given
in~(\ref{postp}) and $\pi(\S^-,\S^+|\theta,\T)$ is defined in
Theorem~\ref{Prop1}(i).
}
\end{thm}

The above characterization of the posterior distribution of $G$
given $\theta$ for models in~(\ref{mixture1}) that is in terms of
two $\S$-paths is less complex than~(or as complex as only when $n,
N-n < 4$) the partition-based characterization~(\ref{postp})~(see
Remark~\ref{remarkbell} for discussion in detail). A proof of the
above two theorems is given in the Appendix.

Given any path ${\S}^+$ and ${\S}^-$ of $n+1$ and $N-n+1$
coordinates, respectively, define
\begin{eqnarray}
  \eta_{0}(\S^+,\S^-) &=&
  \frac{\chi(\mathcal{M}_{1,n}(\S^+),\mathcal{M}_{1,N-n}(\S^-),1)}
  {\chi(\mathcal{M}_{1,n}(\S^+),\mathcal{M}_{1,N-n}(\S^-))}, \label{eta0}\\
\noalign{\noindent for $j=1,\ldots,n$,}
  \eta^+_{j}(\S^+,\S^-) &=&
\frac{\chi(\mathcal{M}_{1,n}(\S^+)\backslash\{S_j^+-S_{j-1}^+\},
  \mathcal{M}_{1,N-n}(\S^-),S_j^+-S_{j-1}^++1)}
{\chi(\mathcal{M}_{1,n}(\S^+),\mathcal{M}_{1,N-n}(\S^-))},
\label{etapos}\\
\noalign{\noindent if $S_j^+>S_{j-1}^+$, 0 otherwise, and, for
$j=1,\ldots,N-n$,}
  \eta^-_{j}(\S^+,\S^-) &=&
  \frac{\chi(\mathcal{M}_{1,n}(\S^+),
  \mathcal{M}_{1,N-n}(\S^-)\backslash\{S_j^--S_{j-1}^-\},S_j^--S_{j-1}^-+1)}
{\chi(\mathcal{M}_{1,n}(\S^+),\mathcal{M}_{1,N-n}(\S^-))}\label{etaneg},
\end{eqnarray}
if $S_j^->S_{j-1}^-$, 0 otherwise.

\begin{cor}\label{cor1}
{\rm Theorems~\ref{Prop1} and~\ref{Prop2} imply that a Bayes
estimate of the unimodal density~(\ref{mixture1}) is given by the
posterior mean given $\theta$ and $\T$,
\begin{equation}\label{EStheta}
    \E[ f( t| G,\theta ) | \theta,\T ] =
    \sum_{\S^+} \sum_{\S^-} a_f(t|\S^-,\S^+,\theta,\T)
     \pi(\S^-,\S^+|\theta,\T)
\end{equation}
where
\begin{eqnarray}
  a_f(t|\S^-,\S^+,\theta,\T)\hspace*{-0.2in}& &= \left[ \eta_{0}(\S^+,\S^-)\, d^+_{\theta,0}(t)
    +\sum_{j=1}^{n} \eta^+_{j}(\S^+,\S^-)\, d_{\theta,j}^+(t|
    \S^+)\right] \I{t > \theta}\hspace*{0.5in}\nonumber\\
   && \quad +\left[ \eta_{0}(\S^+,\S^-)\, d^-_{\theta,0}(t)
    +\sum_{j=1}^{N-n} \eta^-_{j}(\S^+,\S^-)\, d_{\theta,j}^-(t|
    \S^-)\right] \I{t < \theta},\label{af}
\end{eqnarray}
with $d^+_{\theta,0}(t ) = \int_{t-\theta}^{\infty} X^{-1} H(dX)$,
$d_{\theta,0}^-(t ) = \int_{-\infty}^{t-\theta} (-X)^{-1} H(dX)$,
for $j=1,\ldots,n$,
\begin{eqnarray*}
&&d_{\theta,j}^+(t | \S^+) =
\left\{
\begin{array}{lll}
\dfrac{\int_{\max(t-\theta,Y_j)}^{\infty} U^{-(S_j^+-S_{j-1}^++1)}
H(dU)}{\int_{Y_j}^{\infty} U^{-(S_j^+-S_{j-1}^+)} H(dU)}, && S_j^+>S_{j-1}^+,\\
0, && \mbox{otherwise,}
\end{array}
\right. \\
\noalign{\noindent and, for $j=1,\ldots,N-n$,} && d_{\theta,j}^-(t |
\S^-) = \left\{
\begin{array}{lll}
\dfrac{\int_{-\infty}^{\min(t-\theta,Z_j)}
(-V)^{-(S_j^--S_{j-1}^-+1)}
H(dV)}{\int_{-\infty}^{Z_j} (-V)^{-(S_j^--S_{j-1}^-)} H(dV)}, && S_j^->S_{j-1}^-,\\
0, && \mbox{otherwise}.
\end{array}
\right.
\end{eqnarray*}
}
\end{cor}


The above Bayes estimate is a weighted average of the function
$a_f(t|\S^-,\S^+,\theta,\T)$ with respect to
$\pi(\S^-,\S^+|\theta,\T)$. When $H$ is defined by~(\ref{priorH})
and the prior components, $\eta_{0}(\S^+,\S^-)\, d^+_{\theta,0}(t)$
and $\eta_{0}(\S^+,\S^-)\, d^-_{\theta,0}(t)$, vanish, the function
$a_f$ becomes constant between different ordered observations, which
is of the same form as Robertson~(1967)'s maximum likelihood
estimate of the unimodal density when the mode is known as $\theta$.

Dividing the right hand side of~(\ref{product}) by the joint
distribution of $(\VV,\UU,\S^-, \S^+)$ given $\theta$ and $\T$,
given by the last line in~(\ref{postS}), yields the following
analogue of Lemma~2.1 in Ho~(2006a) or Corollary~2.4 in Ho~(2006b)
which states that given $(\S^-,\S^+,\theta,\T)$, $\p$ is uniformly
distributed over all partitions that can be split into $\p^+$ and
$\p^-$ corresponding to the given paths ${\S^+}$ and $\S^-$,
respectively. The above estimator for a unimodal density follows
from the same argument as in Ho~(2006b) to be always less variable
than its counterpart in terms of partitions due to~(\ref{postp}) as
a result of Rao--Blackwellization.

\begin{cor}\label{cor2}
{\rm Consider models in~(\ref{mixture1}). Suppose $\S^+|\theta,\T
\sim \pi^+(\S^+|\theta,\T)$ and\break $\S^-|\S^+,\theta,\T \sim
\pi^-(\S^-|\S^+,\theta,\T)$. Then, there exists a conditional
distribution
$$
    \pi(\p|\S^-,\S^+,\theta,\T) = \frac{1}{|{\mathbb{C}}_{\S^+}|
    |{\mathbb{C}}_{\S^-}|}, \qquad \p = \p^+ \cup \p^-,\p^+ \in
    {\C}_{\S^+}, \p^- \in
    {\C}_{\S^-},
$$
where $|{\mathbb{C}}_{\S}|$ is given by~\mref{card}. }
\end{cor}

Suppose $\theta$ is not known. Theorems~\ref{Prop1} and~\ref{Prop2}
yield the conditional density of $\T$ given $\theta$ to be
proportional to $\sum_{\S^+} \sum_{\S^-} \phi^+_\theta({\S^+},\T)
\,\phi^-_\theta({\S^-},\S^+,\T)$. A standard prior-posterior
updating operation in which the prior distribution of $\theta$ is
$\pi(d\theta)$ results in the following theorem.

\begin{thm}\label{Prop3} {\rm
Assume any prior $\pi(d\theta)$ for $\theta$. Then, the posterior
distribution of $\theta$ given $N$ i.i.d. observations $\T$
from~(\ref{mixture1}) 
is characterized by
\begin{equation}
  \Pr(\theta \in B|\T) = \int_B \sum_{\S^+} \sum_{\S^-}
  \pi(\S^-,\S^+,d\theta|\T)\label{posttheta},
\end{equation}
for any Borel set $B \in \R$, where
\begin{equation}\label{jointAll}
\pi(\S^-,\S^+,d\theta|\T) = \frac{\phi^+_\theta({\S^+},\T)
\,\phi^-_\theta({\S^-},\S^+,\T)\,\pi(d\theta)}{\int_\R\sum_{\S^+}
\sum_{\S^-}\phi^+_\vartheta({\S^+},\T)
\,\phi^-_\vartheta({\S^-},\S^+,\T)\,\pi(d\vartheta)},
\end{equation}
with $\phi^+_\theta({\S^+},\T)$ and $\phi^-_\theta({\S^-},\S^+,\T)$
defined in~(\ref{ZSphi}) and~(\ref{ZSphi2}), respectively. }
\end{thm}

Finally, the posterior distribution of the pair $(G,\theta)$
in~(\ref{mixture1}) can be determined from~(\ref{double}) based on
Theorems~\ref{Prop2} and~\ref{Prop3} and, hence, the posterior
expectation of any functional $h$ of $(G,\theta)$ can be expressible
in terms of a finite sum over two dependent $\S$-paths. In
particular, a Bayes estimate of the unknown unimodal
density~(\ref{mixture1}) follows by letting $h(G,\theta) =
f(t|G,\theta)$ in~(\ref{double}) and applying Corollary~\ref{cor1}.

\begin{thm}\label{Prop4}
{\rm Assume any prior $\pi(d\theta)$ for $\theta$. Then, the
posterior mean of an unimodal density~(\ref{mixture1}) given $N$
i.i.d. observations $\T$ is given by
\begin{equation}\label{Ef}
    \E[ f( t| G,\theta ) |\T ] = \int_\R \sum_{\S^+} \sum_{\S^-}
    a_f(t|\S^-,\S^+,\theta,\T) \pi(\S^-,\S^+,d\theta|\T),
\end{equation}
where $a_f(t|\S^-,\S^+,\theta,\T)$ and $\pi(\S^-,\S^+,d\theta|\T)$
are defined in Corollary~\ref{cor1} and~(\ref{jointAll}),
respectively. }
\end{thm}

\begin{rem}\label{remarkbell}
{\rm As the total number of partitions of $N$ integers, which is the
Bell's exponential number $B_N$, is roughly equal to $N!$, the
complexity of the partition-based characterization~(\ref{postp}),
which relies on the total number of partitions $\p$~(that is, the
number of summands in $\sum_\p$), is roughly equal to $n! \times
(N-n)!$. This quantity is identical to that of Brunner's model, but
it has not been pointed out in Brunner~(1992). Meanwhile, the
complexity of the path-based characterization~(\ref{postS}), which
is based on the number of summands in the double sum
$\sum_{\S^+}\sum_{\S^-}$, depends on $\Lambda_n \times
\Lambda_{N-n}$ where $\Lambda_n$ denotes the total number of
$\S$-paths of $n+1$ coordinates. Hence, its complexity is less than
that of~(\ref{postp}) except when both $n$ and $N-n$ are less than 4
because $\Lambda_n \leq B_n$ for all integers $n$, with equality
only when $n < 4$~(Brunner and Lo~(1989) and Ho~(2002)).
Table~\ref{totalnumber} reveals a ratio between the complexities
of~(\ref{postS}) and~(\ref{postp}) to be as large as $0.02097$ when
$N = 20$. This upper bound on the ratio drops quickly when $N$
increases; for example, the bound is given by $0.00013^2 = 1.69
\times 10^{-8}$ when $N = 40$. }
\end{rem}

\begin{table}[ht]
\renewcommand{\baselinestretch}{1} \tiny \normalsize
\centering\parbox{5in}{\caption {Complexities between path-based and
partition-based characterizations,~(\ref{postS}) and~(\ref{postp}),
versus sample sizes $n$ and $20-n$. }\label{totalnumber} } \vskip
0.2in
\begin{tabular}[h]{rrrcrcr}
\hline\hline $n$& $20-n$ &
\multicolumn{2}{c}{$\Lambda_n \times \Lambda_{20-n}$} 
&
\multicolumn{2}{c}{$B_n \times B_{20-n}$} 
& Ratio in \%  \\ \hline
10 & 10 & 282{,}105{,}616 & & 13{,}450{,}200{,}625 & & 2.097 \\
8 & 12 & 297{,}457{,}160 & & 17{,}444{,}291{,}580 & & 1.705 \\
6 & 14 & 353{,}026{,}080 & & 38{,}752{,}562{,}366 & & 0.911 \\
4 & 16 & 495{,}007{,}380 & & 157{,}202{,}132{,}205 & & 0.315\\
2 & 18 & 1{,}432{,}916{,}100 & & 2{,}046{,}230{,}418{,}477 & &  0.070 \\
0 & 20 & 6{,}564{,}120{,}420 & & 51{,}724{,}158{,}235{,}372 & &  0.013 
\\ \hline\hline
\end{tabular}
\end{table}

\subsection{An illustration with the two-parameter Poisson-Dirichlet
process}\label{sec:PD}

This section illustrates results obtained so far by selecting an
important example of the class of species sampling
models~\mref{randomP}, namely, the two-parameter Poisson-Dirichlet
process~(Pitman and Yor~(1997)). Write the random measure as
$\PY(H;a,b)$ to indicate that its shape probability is $H$ and there
are two shape parameters $0 \leq a <1$ and $b
> -a$. A Dirichlet process with shape measure $\theta H$, $\theta>0$,
corresponds to $\PY(H;0,\theta)$. Selections of $a = \alpha$ and
$b=0$ give a normalized stable law with index $0<\alpha<1$, of which
a simple exponential change of measure gives a normalized
inverse-Gaussian process considered by Lijoi, Mena and
Pr\"unster~(2005). Posterior analysis of models in~(\ref{mixture1})
wherein $G$ is $\PY(H;a,b)$ follow from the previous discussion with
explicit simplifications including $ \ell_{0,k} = (b + N_{k} a)/(b +
k)$ and $\ell_{j,k} = (e_j - a)/(b + k)$ in~(\ref{gBM}),
$$
    \chi(\mathcal{M}_{1,n}(\S)) = 
    \dfrac{\prod_{i=1}^{N(\S)}
    [b + (i-1) a ]
    \prod_{\{j^{\ast}|\S\}} \prod_{i=1}^{S_j - S_{j-1} - 1} (i - a)}
    {\prod_{k=1}^{n} (b + k -1)}
$$
in~(\ref{priorS}--\ref{ZSphi}),
$$
\frac{\chi(\mathcal{M}_{1,n}(\S^+),\mathcal{M}_{1,N-n}(\S^-))}
    {\chi(\mathcal{M}_{1,n}(\S^+))} = \dfrac{\prod_{i=1}^{N(\S^-)}
    [b + (N(\S^+)+i-1) a ]
    \prod_{\{j^{\ast}|\S^-\}} \prod_{i=1}^{S_j^- - S_{j-1}^- - 1} (i - a)}
    {\prod_{k=n+1}^{N} (b + k -1)}
$$
in~(\ref{ZSphi2}), $\eta_{0}(\S^+,\S^-) = [b+ (N(\S^+)+N(\S^-))a] /
(b+N)$, $\eta^+_{j}(\S^+,\S^-) = (S_j^+-S_{j-1}^+-a) / (b+N)$, and
$\eta^-_{j}(\S^+,\S^-) = (S_j^--S_{j-1}^--a) / (b+N)$
in~(\ref{eta0}--\ref{etaneg}), respectively. Last but not least, in
Theorem~\ref{Prop2}, $\CP(dG|\VV,\UU,\S^-, \S^+,\theta,\T)$ is
equivalent to
$$
\sum_{\{j^\ast|\S^+\}}\frac{G^+_j}{G^{\ast}} \delta_{U_j}(\cdot) +
\sum_{\{j^\ast|\S^-\}}\frac{G^-_j}{G^{\ast}} \delta_{V_j}(\cdot) +
\left(1-\sum_{\{j^\ast|\S^+\}}\frac{G^+_j}{G^{\ast}}
-\sum_{\{j^\ast|\S^-\}}\frac{G^-_j}{G^{\ast}}\right)
\CP^{\ast}(\cdot),
$$
where $G^+_j \stackrel{\mbox{\scriptsize ind}}{\sim}
Gamma(S_j^+-S_{j-1}^+-a)$, $G^-_j \stackrel{\mbox{\scriptsize
ind}}{\sim} Gamma(S_j^--S_{j-1}^--a)$, $G^{\ast} =
\sum_{\{j^\ast|\S^+\}}G^+_j + \sum_{\{j^\ast|\S^-\}} G^-_j + G$ with
$G \sim Gamma(b+(N(\S^+)+N(\S^-))a)$, and all variables are mutually
independent of $\CP^{\ast} = \PY(H;a,b+(N(\S^+)+N(\S^-))a)$. Note
that the above new expression of~(\ref{ZSphi2}) is a symmetric
function, yet different from that in~(\ref{ZSphi}), depending on
$\{S_j^- - S_{j-1}^-: S_j^- > S_{j-1}^-, j=1,\ldots,N-n\}$ only.
This allows a straightforward application of the SIP
sampler~(Algorithm~\ref{SIP}) in drawing $\S^-$ given $\S^+$ when
constructing an SIS method~(Algorithm~\ref{bigSIS}) for models
in~(\ref{mixture1}) in the next section.

\section{Sequential Importance Sampling Schemes}\label{sec:SIS}

This section introduces a SIS method~(Kong, Liu and Wong~(1994), Liu
and Chen~(1998) and Liu, Chen and Wong~(1998)) for sampling the
triplets $( \S^-, \S^+,\theta)$ in evaluating/approximating
posterior quantities for models in~(\ref{mixture1}),
like~(\ref{posttheta}) and~(\ref{Ef}), which are expressible in
terms of finite sums of two dependent $\S$-paths. The SIS method is
based on yet another novel SIS~algorithm, called \textit{sequential
importance path~(SIP) sampler}, for sampling one single path at a
time. The SIP sampler is designed in accordance with choosing trial
distributions that mimic the probability kernels for Markov
transitions in the accelerated path~(AP) sampler introduced in
Ho~(2002, 2006a, 2006b) that serves the same purpose.

Generally speaking, the SIP sampler or any other existing SIS method
allows us to draw an $\S$-path of $n+1$ coordinates according to a
probability distribution
\begin{equation}\label{piS} 
\pi(\S)
\propto \phi(\S) 
= |\C_{\S}| \chi(\mathcal{M}_{1,n}(\S)) \prod_{\{j^*|\S\}}
m^{(S_j-S_{j-1})}(Q_j),
\end{equation}
where $|\mathbb{C}_{\S}|$ is defined in~(\ref{card}), $\chi(\cdot)$
is a symmetric function depending on only its arguments, similar
to~(\ref{priorS}), $m^{(S_j-S_{j-1})}(Q_j)$ is a finite real-valued
function depending on $S_j-S_{j-1}$ and $Q_j$ only, and $Q_1,
\ldots, Q_n$ is a decreasing/increasing sequence in $\R$. An
inefficient SIS method proposed by Ho~(2002, Section~4) consists of
$n-1$ recursive determinations of one coordinate of the path $\S$ at
a time in an ascending order conditioning on all previously
determined coordinates according to a trial distribution
\begin{equation}\label{naivet}
    \Pr(S_r = s_r|S_1 = s_1,\ldots,S_{r-1} = s_{r-1})
    := t_r(s_r|s_1,\ldots,s_{r-1}) \propto
    \phi(\s_{r,s_r}^{\ast})
\end{equation}
for $r=1,\ldots,n-1$, where $\s_{r,s_r}^{\ast} =
(0,s_1,\ldots,s_{r-1},s_r,r+1)$ is a path of $r+2$ coordinates.
After step $n-1$, a path $\s = (0,s_1,s_2,\ldots,s_{n-1},n)$ drawn
with probability $t_{n-1}(\s) = \prod_{r=1}^{n-1}
t_r(s_r|s_1,\ldots,s_{r-1})$ can then be treated as a Monte Carlo
sample from~(\ref{piS}) after being properly weighted by an
importance sampling weight $\phi(\s) / t_{n-1}(\s)$. However, it
turns out that the above scheme is practically not efficient in
evaluating sums over $\S$-paths. In general, this is directly
related to the discrepancy between the trial distribution
$t_i(\cdot|\cdot)$ in~(\ref{naivet}) and the true conditional
distribution of $S_r$ given $S_1,\ldots,S_{r-1}$ derivable from the
target distribution $\pi(\S)$~(Liu and Chen~(1998)). Noticing that
the ``transition'' in~(\ref{naivet}) is equivalent to determinations
of the two increments, $S_r - S_{r-1}$ and $r+1 - S_r$, of the path
at locations $r$ and $r+1$, respectively, our idea is to replace
location $r+1$ by some other latter location $q$, which parallels
the idea of constructing the AP sampler adopted in Ho~(2002,
2006a,b) when improving on an inefficient Gibbs chain. Let $I_{0}=0$
and $I_{n}=n$ and denote $\{I_{1},\ldots,I_{n-1}\}$ as a random
permutation of the integers $\{1,2,\ldots,n-1\}$, such that $D_r =
\{I_0\} \cup \{I_1,\ldots,I_r\} \cup \{I_n\}$ consists of all
determined coordinates of the $\S$-path after step $r$ of the SIP
sampler, for $r=1,\ldots,n-1$,

\begin{algorithm}[Sequential importance path~(SIP) sampler]\label{SIP}
{\rm An efficient SIS method for sampling an $\S$-path of $n+1$
coordinates from $\pi(\S)$ given in~(\ref{piS}), the SIP sampler,
consists of recursive applications of the following SIS steps for
$r=1,\ldots,n-1$:
\begin{enumerate}
    \item[A.] Given $D_{r-1}$, let
    $p=\max\{I_{j} \in D_{r-1}:I_{j}<I_r\}$ and
$q=\min\{I_{j}\in D_{r-1}:I_{j}>I_r\}$. Determine $S_{I_r} = k$, for
$k = S_p, S_p+1, \ldots, \min(I_r,S_q)$, according to a distribution
    \begin{equation}
    \kappa_{r}(k|\{S_h:h \in D_{r-1}\})
    \propto \phi(\S_{I_r,k}^{\ast}),
    \end{equation}
    where $\S_{I_r,k}^{\ast} = (0,S^{\ast}_1,\ldots,S^{\ast}_{I_r-1},S^{\ast}_{I_r},
    S^{\ast}_{I_r+1},\ldots,S^{\ast}_{n-1},n)$ is a path of $n+1$ coordinates
    such that $S_{I_r}^{\ast} = k$ and for $i=1,\ldots,I_{r}-1,I_{r}+1,\ldots,n-1$,
    $S_i^{\ast} = S_{I_h}$ if $i = I_h \in D_{r-1}$; otherwise, $S_i^{\ast}
    = S_{i-1}^{\ast}$~(see Remark~\ref{prob} for explicit
    expressions of $\kappa_{r}(k|\{S_h:h \in D_{r-1}\})$ for
    different values of $k$).
    \item[B.] Compute $\kappa_{r}(k|\{S_h:h \in D_{r-1}\})$, equals
    $\phi(\S^\ast_{I_r,k})$
    multiplied by the appropriate constant of proportionality, for the chosen
    value $k$ of $S_{I_r}$.
\end{enumerate}
}
\end{algorithm}

After step $n-1$, we obtain a random path $\S =
(0,S_1,S_2,\ldots,S_{n-1},n)$ distributed as the trial distribution
\begin{equation}\label{kappa}
\kappa_{n-1}(\S) = \prod_{r=1}^{n-1} \kappa_{r}(S_{I_r}|\{S_h:h \in
D_{r-1}\}).
\end{equation}
Hence, its importance sampling weight is given by $w_{n-1}(\S) =
\phi(\S) / \kappa_{n-1}(\S)$. Given $M$ i.i.d. draws,
$\S_{(1)},\ldots,\S_{(M)}$ with respective importance sampling
weights \break $w_{n-1}(\S_{(1)}), \ldots,w_{n-1}(\S_{(M)})$, from
the SIP sampler based on different permutations
$\{I_{1},\ldots,I_{n-1}\}$ of the $n-1$ integers, any sum over
$\S$-paths/expectation of any functional $h(\S)$ with respect to the
probability distribution $\pi(\S)$, $\eta_h = \sum_{\S} h(\S)
\pi(\S)$, can be approximated by
\begin{equation}\label{SIS1}
\eta_h^{M} = \frac{\sum_{i=1}^M h(\S_{(i)})\,
w_{n-1}(\S_{(i)})}{\sum_{i=1}^M w_{n-1}(\S_{(i)})}.
\end{equation}

\begin{rem}{\rm We remark that there are two major differences between the SIP
sampler and the inefficient SIS method which intuitively explain why
the SIP sampler is more efficient. On one hand, the coordinates of
the path are determined in a random order in the SIP sampler, but
not in an ascending order or any other pre-determined order. This
arrangement is desired and crucial, as it results in determination
of an increment at a location possibly latter than $r+1$ in step
$r$, which is the idea behind the success of the AP sampler. On the
other hand, each trial distribution $\kappa_{r}(S_{I_r}|\{S_h:h \in
D_{r-1}\})$ in the SIP sampler mimics the transition probabilities
in the AP sampler, in the sense that it is proportional to the
probability of a path of $n+1$ coordinates for any $r =
1,\ldots,n-1$, rather than the probability of a path of number of
coordinates varying with $r$. }\end{rem}

\begin{rem}\label{prob}{\rm In the SIP sampler~(Algorithm~\ref{SIP}), the trial distribution
$\kappa_{r}(k|\{S_h:h \in D_{r-1}\})$ is explicitly proportional
to\vspace*{-0.05in}
\begin{eqnarray}
    && \frac{I_r-S_p}{S_q-S_p-1}
    \, \chi(\mathcal{M}_{1,p}(\S_{I_r,S_p}^{\ast}),
    S_q-S_p, \mathcal{M}_{q+1,n}(\S_{I_r,S_p}^{\ast})) \, m^{(S_q - S_p)}(Q_q) \nonumber \\
\noalign{\noindent if $k = S_p$, or}
    &&    \I{I_r \geq S_q}\frac{I_r-S_p}{S_q-S_p-1}
    \prod_{i=S_p+1}^{S_q-1} \left(\frac{I_r - i}{q-i}\right) \nonumber \\
    && \qquad \times \chi(\mathcal{M}_{1,p}(\S_{I_r,S_q}^{\ast}),
    S_q-S_p, \mathcal{M}_{q+1,n}(\S_{I_r,S_q}^{\ast}))\, m^{(S_q - S_p)}(Q_p) \nonumber\\
\noalign{\noindent if $k = S_q$, or}
    &&    \hspace*{-0.2in}\binom{S_q-S_p-2}{S_q-k-1}
    \prod_{i=I_r+1}^{q-1} \left(\frac{i-k}{I_r-S_p}\right)\nonumber\\
    && \times \chi(\mathcal{M}_{1,p}(\S_{I_r,k}^{\ast}),
    k-S_p,S_q-k, \mathcal{M}_{q+1,n}(\S_{I_r,k}^{\ast}))\,
    m^{(k-S_p)}(Q_p)\,m^{(S_q - k)}(Q_q).\nonumber
\end{eqnarray}
if $k=S_p+1,\ldots,\min(I_r,S_q-1)$.
}\end{rem}


\begin{algorithm}\label{bigSIS}
{\rm An SIS method that samples $( \S^-, \S^+,\theta)$
from~(\ref{jointAll}) consists of three major steps:
\begin{enumerate}
    \item[(i)] Sample $\theta$ according to a density $\rho(\theta) >0$, $\theta \in \R$. Then, define $\Y$ and $\Z$
    accordingly based on~(\ref{obs}). Also, choose random permutations of
    $\{1,\ldots,n-1\}$ and $\{1,\ldots,N-n-1\}$.
    \item[(ii)] Given $\theta$, determine $\S^+$ by applying Algorithm~\ref{SIP}
    with function $\phi(\S)$ defined by $\phi^+_\theta(\S^+,\T)$ in~(\ref{ZSphi}).
    Obtain $\kappa_{n-1}(\S^+|\theta)$ according to~(\ref{kappa}).
    \item[(iii)] Given $(\S^+,\theta)$, determine $\S^-$ by applying Algorithm~\ref{SIP}
    with function $\phi(\S)$ defined by
    $\phi^-_\theta(\S^-,\S^+,\T)$ in~(\ref{ZSphi2}), provided that the ratio
    $\chi(\mathcal{M}_{1,n}(\S^+),\mathcal{M}_{1,N-n}(\S^-)) /
    \chi(\mathcal{M}_{1,n}(\S^+))$ is a
    symmetric function of $\mathcal{M}_{1,N-n}(\S^-)$.
    Obtain $\kappa_{N-n-1}(\S^-|\S^+,\theta)$ according to~(\ref{kappa}).
\end{enumerate}
}
\end{algorithm}

After a total of $N-1$ sub-steps, we obtain a random sample of
$(\S^-,\S^+,\theta)$ distributed as the trial distribution
$
\kappa_{N-n-1}(\S^-|\S^+,\theta)\times \kappa_{n-1}(\S^+|\theta)
\times \rho(\theta) .
$ 
If $\pi(d\theta) = \pi(\theta)d\theta$, its importance sampling
weight is given by
$$
w_{N-1}(\S^-,\S^+,\theta) = \frac{\phi^+_\theta({\S^+},\T)
\,\phi^-_\theta({\S^-},\S^+,\T) \,\pi(\theta)}
{\kappa_{N-n-1}(\S^-|\S^+,\theta)\, \kappa_{n-1}(\S^+|\theta) \,
\rho(\theta)}.
$$
We remark that it is possible, indeed more desired in terms of
efficiency of the SIS method, that the sequence in sampling the two
paths in steps~(ii) and~(iii) can be randomized based on
appropriate, but slight, modifications of the function $\phi(\S)$ in
applying Algorithm~\ref{SIP}. That is, there is one-half probability
that $\S^+$ is sampled before $\S^-$ as stated in
Algorithm~\ref{bigSIS}; otherwise, $\S^-$ is sampled before $\S^+$.

\begin{cor}\label{bigapprox}{\rm
Posterior quantities for models in~(\ref{mixture1}),
like~(\ref{posttheta}) and~(\ref{Ef}), which are expressible as
$$
\gamma_h =  \int_\R \sum_{\S^-}\sum_{\S^+} h(\S^-,\S^+,\theta)\,
\pi(\S^-,\S^+,d\theta|\T)
$$
can be approximated by
\begin{equation}\label{Efestimate}
\gamma_h^{M} = \frac{\sum_{i=1}^M
h(\S^-_{(i)},\S^+_{(i)},\theta_{(i)})\,
w_{N-1}(\S^-_{(i)},\S^+_{(i)},\theta_{(i)})}{\sum_{i=1}^M
w_{N-1}(\S^-_{(i)},\S^+_{(i)},\theta_{(i)})},
\end{equation}
where $(\S^-_{(1)},\S^+_{(1)},\theta_{(1)}),\ldots,
(\S^-_{(M)},\S^+_{(M)},\theta_{(M)})$ is a sequence of $M$ i.i.d.
samples from~(\ref{jointAll}) with respective importance sampling
weights $w_{N-1}(\S^-_{(1)},\S^+_{(1)},\theta_{(1)}), \ldots,
w_{N-1}(\S^-_{(M)},\S^+_{(M)},\theta_{(M)})$, obtained by carrying
out Algorithm~\ref{bigSIS} independently for a large number of times
$M$. }
\end{cor}

\section{Numerical Results}\label{sec:sim}

This section concerns practical applications of our methodology. For
purpose of illustration, $G$ is selected to be the two-parameter
Poisson-Dirichlet process as the corresponding results are discussed
in Section~\ref{sec:PD}. The idea of conjugacy suggests $H(\cdot)$
of the measure $\PY(H;a,b)$ to be related to a Pareto distribution.
In particular, we chose the following mixture of two Pareto random
variables, symmetrical about zero, that is,
\begin{equation}\label{priorH}
H(dX) = \frac{\alpha \delta^{\alpha}}{2(-X)^{\alpha+1}}
\I{X<-\delta} dX + \frac{\alpha \delta^{\alpha}}{2X^{\alpha+1}}
\I{X>\delta} dX, \qquad \alpha,\delta>0,
\end{equation}
such that it results in
$$
\int_Y^\infty X^{-\nu} H(dX) = \int_{-\infty}^{-Y} (-X)^{-\nu} H(dX)
= \frac{\alpha \delta^{\alpha}}
{2(\alpha+\nu)\max(|Y|,\delta)^{\alpha+\nu}}.
$$
for $Y > 0$ and any positive integer $\nu$, which are necessary in
implementation of Algorithm~\ref{bigSIS}~(or Algorithm~\ref{SIP}).
For purpose of ``deflating'' the prior belief, we choose $\alpha =
\delta = 0.000001$. Due to the same reason, the prior $\pi(d\theta)$
is chosen to be uniformly distributed on a reasonably large interval
on $\R$ such that all observations are included. The sequence in
which the coordinates of the $\S$-paths are determined, say,
$\{I_1,\ldots,I_{n-1}\}$ for a path of $n+1$ coordinates, is
randomized in every application of the sequential algorithms.
Likewise, the determinations of the two paths are also randomized in
Algorithm~\ref{bigSIS}. Last but not least, the Monte Carlo size $M
= 1000$.

\subsection{Resolution of the SIP sampler}\label{sec:SIPresolution}
This section addresses the performance of the SIP sampler which
directly affects the SIS method~(Algorithm~\ref{bigSIS}) for
estimating a unimodal density. Based on a fixed and known mode
$\theta_0$, our interest is to estimate the unimodal
density~(\ref{mixture1}) with $G$ taken to be the two-parameter
Poisson-Dirichlet process with $H$ in~(\ref{priorH}), $a = 0$ and $b
= 1$, given as in~(\ref{EStheta}) with $a_f(t|\S^-,\S^+,\theta,\T)$
defined by $\theta = \theta_0$ together with the simplifications
discussed in Section~\ref{sec:PD}. To approximate the posterior
mean, steps~(ii) and~(iii) in Algorithm~\ref{bigSIS}\footnote{As
discussed after the introduction of Algorithm~\ref{bigSIS}, the
sequence of determinations of the two paths -- $\S^+$ first or
$\S^-$ first -- is randomized to achieve a higher efficiency.},
which are essentially two sequential applications of the SIP
sampler, are implemented based on the known mode $\theta_0$. In
particular, the convergence property of the approximated density
estimate as the sample size $N$ increases is studied.

Based on nested samples of sizes $N=500$, $1000$ and $3000$ from a
unimodal density with $[-1,0]$ as modal interval~(Wegman~1970a)
given by
$$
\lambda_1(t) =\left\{
\begin{array}{lcl}
0.02 &  & -7< t\leq -2 \\
0.1 &  & -2< t\leq -1\\
0.4 &  & -1< t \leq 0\\
0.4 \exp(-t) & & t > 0\\
0 & & \mbox{otherwise,}
\end{array}
\right.
$$
weighted averages that approximate the
posterior mean of the unimodal density conditioning on $\theta =
\theta_0$, given as in~(\ref{Efestimate}),
$$
\gamma_{a_f}^{M}(t|\theta_0) = \frac{\sum_{i=1}^M
a_f(t|\S^-_{(i)},\S^+_{(i)},\theta_0,\T)\,
w_{N-2}(\S^-_{(i)},\S^+_{(i)}|\theta_0)}{\sum_{i=1}^M
w_{N-2}(\S^-_{(i)},\S^+_{(i)}|\theta_0)},
$$
where $w_{N-2}(\S^-_{(i)},\S^+_{(i)}|\theta_0)$ is the importance
sampling weight of the pair $(\S^-_{(i)},\S^+_{(i)})$ resulted from
steps~(ii) and~(iii) of Algorithm~\ref{bigSIS}, are displayed at the
left columns in Figures~\ref{fig:SIP1}-\ref{fig:SIP3} for $\theta_0
= -1, -0.5~(\mbox{center mode})$, and $0$, respectively. The whole
procedure is repeated for the two-parameter Poisson-Dirichlet
process with $a=0.9$ and $b=100$. The density estimates based on the
three selected values of the mode are depicted at the right columns
in Figures~\ref{fig:SIP1}-\ref{fig:SIP3}. The graphs echo the fact
that the approximated Bayes estimate of the unimodal density,
$\gamma_{a_f}^{M}(t|\theta_0)$, tends to the ``true'' unimodal
density $\lambda_1(t)$ as sample size increases~(from top to bottom
in the figures) regardless of the two sets of parameters for
$G$~(between columns in the figures). When $N$ is large, there is
not much difference among density estimates based on different
modes.

\subsection{Resolution of the SIS method~(Algorithm~\ref{bigSIS})}\label{sec:bigSISresolution}
The practicality of the SIS method~(Algorithm~\ref{bigSIS}) for
estimation of a unimodal density and its mode is addressed in this
section. To estimate the unimodal density~(\ref{mixture1}) with $G$
taken to the two-parameter Poisson-Dirichlet process with $a=0$ and
$b=1$, Algorithm~\ref{bigSIS} based on $\rho(\theta)$ as a standard
normal density is implemented independently for $M=1000$ number of
times to produce random samples of $(\S^-,\S^+,\theta)$ with
importance sampling weight $w_{N-1}(\S^-,\S^+,\theta)$. According to
Corollary~\ref{bigapprox}, two weighted averages, defined as
in~(\ref{Efestimate}),
$$
\theta^{M} = \frac{\sum_{i=1}^M \theta_{(i)}\,
w_{N-1}(\S^-_{(i)},\S^+_{(i)},\theta_{(i)})}{\sum_{i=1}^M
w_{N-1}(\S^-_{(i)},\S^+_{(i)},\theta_{(i)})}
$$
and
$$
\gamma_{a_f}^{M}(t) = \frac{\sum_{i=1}^M
a_f(t|\S^-_{(i)},\S^+_{(i)},\theta_{(i)},\T)\,
w_{N-1}(\S^-_{(i)},\S^+_{(i)},\theta_{(i)})}{\sum_{i=1}^M
w_{N-1}(\S^-_{(i)},\S^+_{(i)},\theta_{(i)})},
$$
are used to approximate the Bayes estimates~(the posterior mean
given $N$ observations) of the unknown mode $\theta$ and the unknown
unimodal density, respectively.

The unimodal density $\lambda_1(t)$ in the previous section and
another two unimodal densities are
chosen as test densities. 
They are,
$$
\lambda_2(t) = \left\{
\begin{array}{lcl}
0.02 &  & -7< t\leq -2 \\
0.25 &  & -2< t\leq 0\\
0.5 &  & 0 < t \leq 0.1\\
0.1 & & 0.1 < t \leq 2.5\\
0 & & \mbox{otherwise,}
\end{array}
\right.
$$
and
$$
\lambda_3(t) = \frac{12}{13}\left[  \zeta\left(1.5x\right)
\I{-\infty<x<0}+\zeta\left( \frac{x}{1.5}\right) \I{ 0<x<\infty}
\right],
$$
where $\zeta( \cdot)$ is the density function of a standard Cauchy
random variable. These three densities behave quite differently from
one another in the sense that they have modal interval of length 1,
modal interval of shorter length 0.1, and a unique mode at zero,
respectively.

Density estimates $\gamma_{a_f}^{M}(t)$ based on nested samples of
sizes $N=500$, $1000$ and $2000$ from the three unimodal densities
are depicted in the left columns of
Figures~\ref{fig:bigSIS1}-\ref{fig:bigSIS3}, respectively, while
mode estimates $\theta^M$ are presented in
Table~\ref{thetaestimate}. The whole procedure is repeated with
$\rho(\theta)$ as a less diffuse normal density with mean 0 and
standard deviation 1/4. The resulting density estimates are depicted
in the right columns of Figures~\ref{fig:bigSIS1}-\ref{fig:bigSIS3},
while mode estimates are appended in Table~\ref{thetaestimate}. It
is evident from the mode estimates in Table~\ref{thetaestimate},
especially when $N$ is not large, that approximation results based
on $\rho(\theta)$ with a smaller standard deviation are much better
than the others. This is also supported by
Figures~\ref{fig:bigSIS1}-\ref{fig:bigSIS3}; for instance, the peak
at the modal interval $[0,0.1]$ is not well-captured even when $N =
2000$~(graph at the bottom-left in Figure~\ref{fig:bigSIS2}). This
phenomenon can be addressed by the well-known fact~(Kong, Liu and
Wong~(1994) and Liu, Chen and Wong~(1998)) that efficiency of any
SIS method depends heavily on whether the initial trial
distributions in its early steps/stages is close to the true
conditional distributions. Hence, a good choice of $\rho(\theta)$ in
step~(i) of Algorithm~\ref{bigSIS} directly affects the efficiency
of the SIS method.

\begin{table}[ht]
\renewcommand{\baselinestretch}{1.3} \tiny \normalsize
\centering\parbox{5in}{\caption {Weighted average estimates of the
mode}\label{thetaestimate} }
\vskip 0.2in
\begin{tabular}{ccccc}
  \hline\hline
  && \multicolumn{3}{c}{Unimodal Density} \\ \cline{3-5}
  $\rho(\theta)$ & $N$ & $\lambda_1(t)$ & $\lambda_2(t)$ & $\lambda_3(t)$ \\ \hline
  & 500 & -1.249450 & -1.538695 & -0.230042 \\
  $N(0,1)$& 1000 & 1.068161 & 0.079308 & 0.815681 \\
  & 2000 & -0.999335 & 0.037052 & -0.615350 \\ \hline
  & 500 & 0.165998 & 0.138150 & -0.071292 \\
  $N(0,0.25^2)$& 1000 & 0.199645 & 0.143027 & 0.013269 \\
  & 2000 & 0.101668 & -0.071598 & -0.271294\\ \hline\hline
  \multicolumn{2}{c}{True mode} & $[-1,0]$ & $[0,0.1]$ & 0\\
  \hline\hline
\end{tabular}
\end{table}

\begin{figure}[t]
\centering\includegraphics[width=6in]{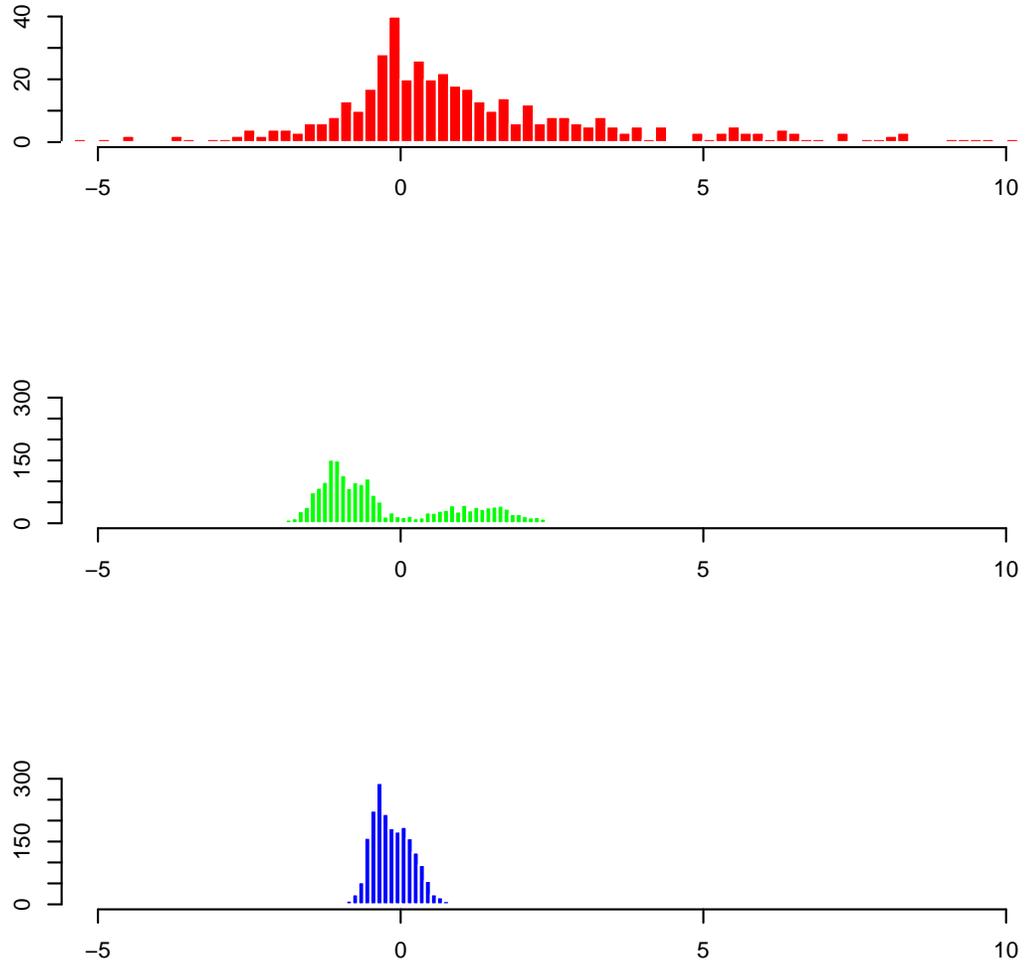}\\
\renewcommand{\baselinestretch}{1}
  \parbox{5in}
  {\caption{\quad Histogram of 500 observations simulated from
  $\lambda_3(t)$ and histograms of the resulting Bayes
  estimates of $\theta$ by 2000 replications of
  Algorithm~\ref{bigSIS} based on a $N(0,1)$
  density and a
  $N(0,0.25^2)$ density for $\rho(\theta)$~(from top to bottom).}
  \label{fig:histtheta}}
\end{figure}

\begin{figure}[t]
\centering\includegraphics[width=6in]{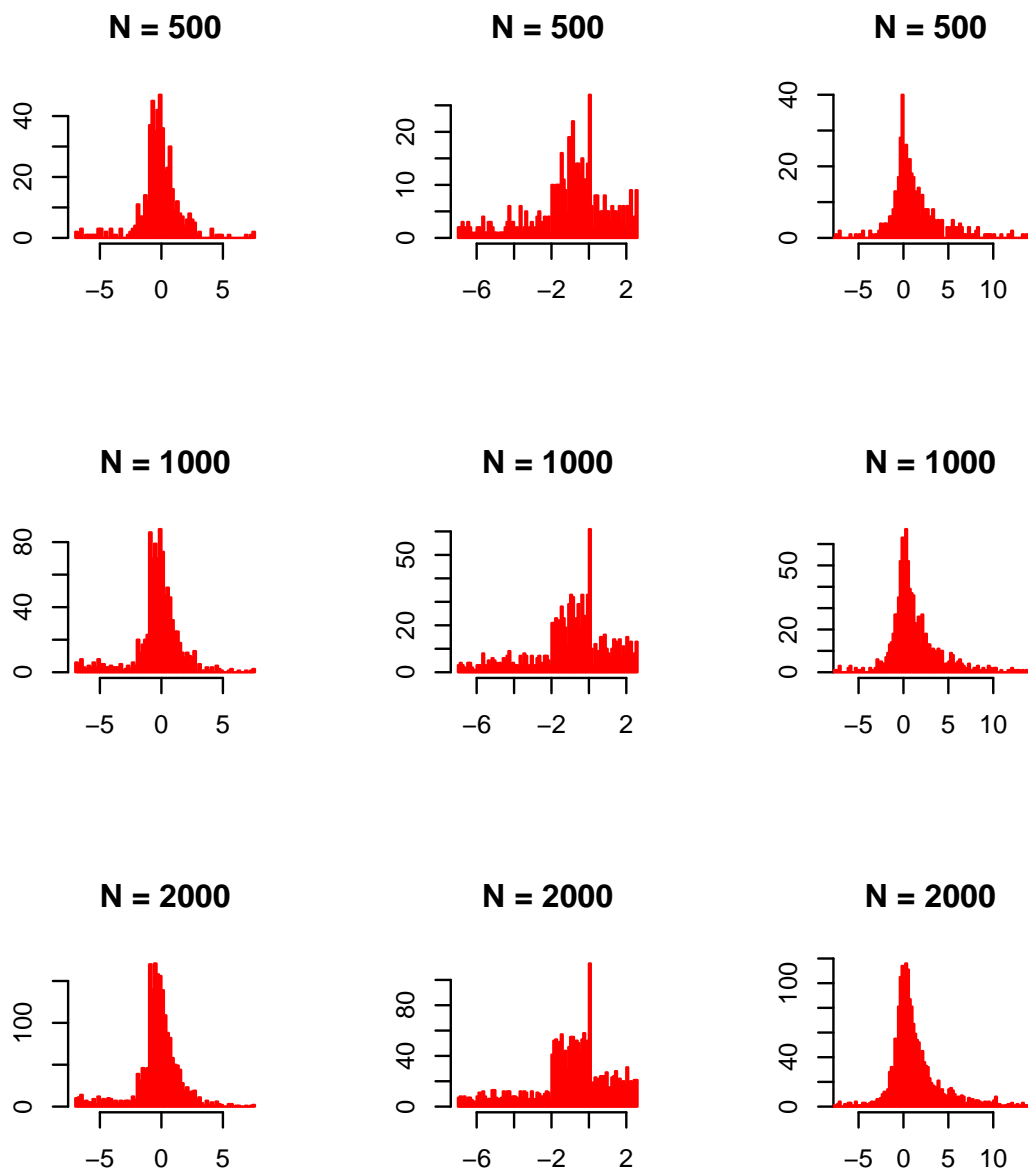}\\
\renewcommand{\baselinestretch}{1}
  \parbox{5in}
  {\caption{\quad Histograms of the data simulated from unimodal densities
  $\lambda_1(t)$~(left column), $\lambda_2(t)$~(middle column), and
  $\lambda_3(t)$~(right column).}
  \label{fig:hist}}
\end{figure}

To explore the selection issue of $\rho(\theta)$, we carry out a
large-sample study by replicating the above procedure to estimate
the mode of the unimodal density $\lambda_3(t)$ based on $N=500$
observations. Histograms of the 2000 independent Bayes estimates of
$\theta$ based on the two different $\rho(\theta)$'s are plotted in
Figure~\ref{fig:histtheta}. It is clear from the graph in the last
row based on a standard normal density for $\rho(\theta)$ does not
give convincing posterior estimates of the mode. On the contrary,
the graph in the second row shows that the true mode is
well-captured when $\rho(\theta)$ is less diffuse. This deficiency
can be understood by looking at the histogram of the 500
observations in the first row in Figure~\ref{fig:histtheta}; the
posterior distribution of the mode $\theta$ should be quite
concentrated around zero and, hence, the choice of a standard normal
density for $\rho(\theta)$ may be far too diffuse. Note that
regarding the results for the first two test densities, the less
diffuse normal density with standard deviation 1/4, symmetrical
about zero, is not really close to the posterior distribution of
$\theta$ at all based on the histograms of the data in
Figure~\ref{fig:hist}. This implies that it is not necessary to set
$\rho(\theta)$ to be extremely close to the true posterior
distribution of $\theta$, which is characterized in
Theorem~\ref{Prop3}. In conclusion, we suggest setting
$\rho(\theta)$ to be a density which is not too diffuse around the
mode~(based on information from the histogram of the data) in
applying the SIS method~(Algorithm~\ref{bigSIS}) for estimating
unimodal densities.

\section*{Appendix: Proof of Theorem~\ref{Prop1} and~\ref{Prop2}}

\begin{proof}
{ Suppose $\theta$ is given. Theorem~2 in Ishwaran and James~(2003)
states that the law of $G$ in~(\ref{mixture1}) given $\theta$ and
$N$ i.i.d. observations $\T$ is characterized by
\begin{eqnarray}
  &&\hspace*{-0.5in}\int_{\CM}g(G)\CP(dG|\theta,\T) \nonumber\\
  &&\hspace*{-0.25in}= \sum_{\p} \left[
  \int_{\R^{\Np}} \left\{\int_\CM g(G) \CP(dG|\X^{\ast},\p,\theta,\T)\right\}
  \prod_{k=1}^{\Np} \mu(dX^{\ast}_k|C_k) \right] \pi(\p|\theta,\T)\label{postp}
\end{eqnarray}
for any nonnegative or integrable function $g$, wherein
$\CP(dG|\X^{\ast},\p,\theta,\T)$ is determined by Lemma~1 in
Ishwaran and James~(2003), and $ \prod_{k=1}^{\Np}
\mu(dX^{\ast}_k|C_k) \pi(\p|\theta,\T)$ is equivalent in
distribution to the posterior distribution of $\X$ given $\theta$ as
discussed in Theorem~1 in Ishwaran and James~(2003), where, for
$k=1,\ldots,\Np$, $\mu(dX^{\ast}_k|C_k)$ is proportional to
\begin{equation}\label{jointp}
\varphi_k(dX_k^{\ast})= \frac{1}{(X_k^{\ast})^{e_k}}
    \left[\I{0< \max_{j \in C_k} T_j-\theta \leq X_k^{\ast}} -
    \I{X_k^{\ast} \leq \min_{j \in C_k} T_j-\theta <0} \right]
    H(dX_k^{\ast})
\end{equation}
and $ \pi(\p|\theta,\T) = \pi(\p)\prod_{k=1}^{\Np}
\int\varphi_k(dX^{\ast}_k) / \sum_\p \pi(\p)\prod_{k=1}^{\Np} \int
\varphi_k(dX^{\ast}_k)$.

Splitting $\p$ into $\p^{+}$ and $\p^{-}$ as discussed in
Section~\ref{sec:post} and re-expressing the kernels according
to~(\ref{twocases}) yield
\begin{eqnarray}\label{varphi}
\varphi_k(dX_k^{\ast}) = \left\{
\begin{array}{lll}
    \varphi^{+}_k(dX_k^{\ast}) = (X_k^{\ast})^{-e_k}
    \I{0< \displaystyle\max_{j \in C_k} Y_j \leq X_k^{\ast}}
    H(dX_k^{\ast}) && k \leq N(\p^+)\\
    \varphi^{-}_k(dX_k^{\ast}) = (-X_k^{\ast})^{-e_k}
    \I{X_k^{\ast} \leq \displaystyle\min_{j \in C_k} Z_j <0} H(dX_k^{\ast})& & k >
    N(\p^+)
\end{array}
\right.
\end{eqnarray}
and
\begin{eqnarray}
\pi(\p|\theta,\T) && \hspace*{-0.25in}=
\frac{\pi(\p)\left[\prod_{i=1}^{N(\p^+)}
\int\varphi_i^{+}(dX_i^{\ast})\right]\left[\prod_{j=N(\p^+)+1}^{\Np}
\int\varphi_j^{-}(dX_j^{\ast})\right]}
    {\sum_\p \pi(\p)\left[\prod_{i=1}^{N(\p^+)}
\int\varphi^{+}_i(dX_i^{\ast})\right]\left[\prod_{j=N(\p^+)+1}^{\Np}
\int\varphi^{-}_j(dX_j^{\ast})\right]}\nonumber \\
&& \hspace*{-0.25in}= \frac{\psi^-(\p^-|\p^+,\theta,\T) \times
\psi^+(\p^+|\theta,\T)}{\sum_{\p^+}
\left\{\sum_{\p^-}\psi^-(\p^-|\p^+,\theta,\T)
\right\}\psi^+(\p^+|\theta,\T)}\label{pip2},
\end{eqnarray}
where $\psi^+(\p^+|\theta,\T) = \pi(\p^+)
\prod_{i=1}^{N(\p^+)}\int\varphi^{+}_i(dX_i^{\ast})$ defines a
posterior distribution of $\p^+$ of the $n$ positive observations
$\Y$ given $\theta$ and $\psi^-(\p^-|\p^+,\theta,\T) =
\pi(\p^-|\p^+)\prod_{j=N(\p^+)+1}^{\Np}
\int\varphi^{-}_j(dX_j^{\ast})$ defines a~(conditional) posterior
distribution of $\p^-$ of the remaining $N-n$ negative observations
given $(\p^+$, $\theta)$. The equality in~(\ref{pip2}) follows
from~(\ref{symmetric2}) and re-writing $\sum_\p = \sum_{\p^+}
\sum_{\p^-}$. Then, combining~(\ref{varphi}) and~(\ref{pip2}) gives
\begin{eqnarray}
\prod_{k=1}^{\Np} \mu(dX^{\ast}_k|C_k) \pi(\p|\theta,\T) &\propto&
    \left[\prod_{j=N(\p^+)+1}^{\Np}
    \varphi^{-}_j(dX_j^{\ast})
    \psi^-(\p^-|\p^+,\theta,\T)\right]\nonumber\\
    && \qquad \qquad \times
    \left[ \prod_{i=1}^{N(\p^+)}
    \varphi^{+}_i(dX_i^{\ast})
    \psi^+(\p^+|\theta,\T)\right].\label{product}
\end{eqnarray}
Theorem~2.1 in Ho~(2006b) yields that the law of
$X^{\ast}_1,\ldots,X^{\ast}_{N(\p^+)},\p^+|\theta,\T$~(proportional
to the last term above) is equivalent to the law of $\UU,
\S^+|\theta,\T$ defined by~(\ref{ZSphi}) and~(\ref{alphaSj}).
Utilizing the symmetric properties of $\pi(\p^-|\p^+)$
in~(\ref{chistar}) and~(\ref{priorS}) and applying Theorem~2.1 in
Ho~(2006b) yield the equivalence between the law of
$X^{\ast}_{N(\p^+)+1},\ldots,X^{\ast}_{\Np},\p^-|\p^+,\theta,\T$,
proportional to the first term at the right hand side
of~(\ref{product}), and the law of $\VV, \S^-|\S^+,\theta,\T$
defined by~(\ref{ZSphi2}) and~(\ref{alphaSj2}), completing the proof
of Theorem~\ref{Prop1}. The result in Theorem~\ref{Prop2} follows as
a result of Theorem~\ref{Prop1} by recognizing the equality in
distribution between $\CP(dG|\X^{\ast},\p,\theta,\T)$ and
$\CP(dG|\VV,\UU,\S^-, \S^+,\theta,\T)$ in~(\ref{postp}).}
\end{proof}






\vskip0.5in {\smc \Tabular{ll}
Man-Wai Ho\\
Department of Statistics and Applied Probability \\
National University of Singapore\\
6 Science Drive 2\\
Singapore 117546\\
Republic of Singapore\\
\rm E-mail: stahmw\at nus.edu.sg\\
\EndTabular }

\begin{figure}[h]
  \centering\includegraphics[width=6in]{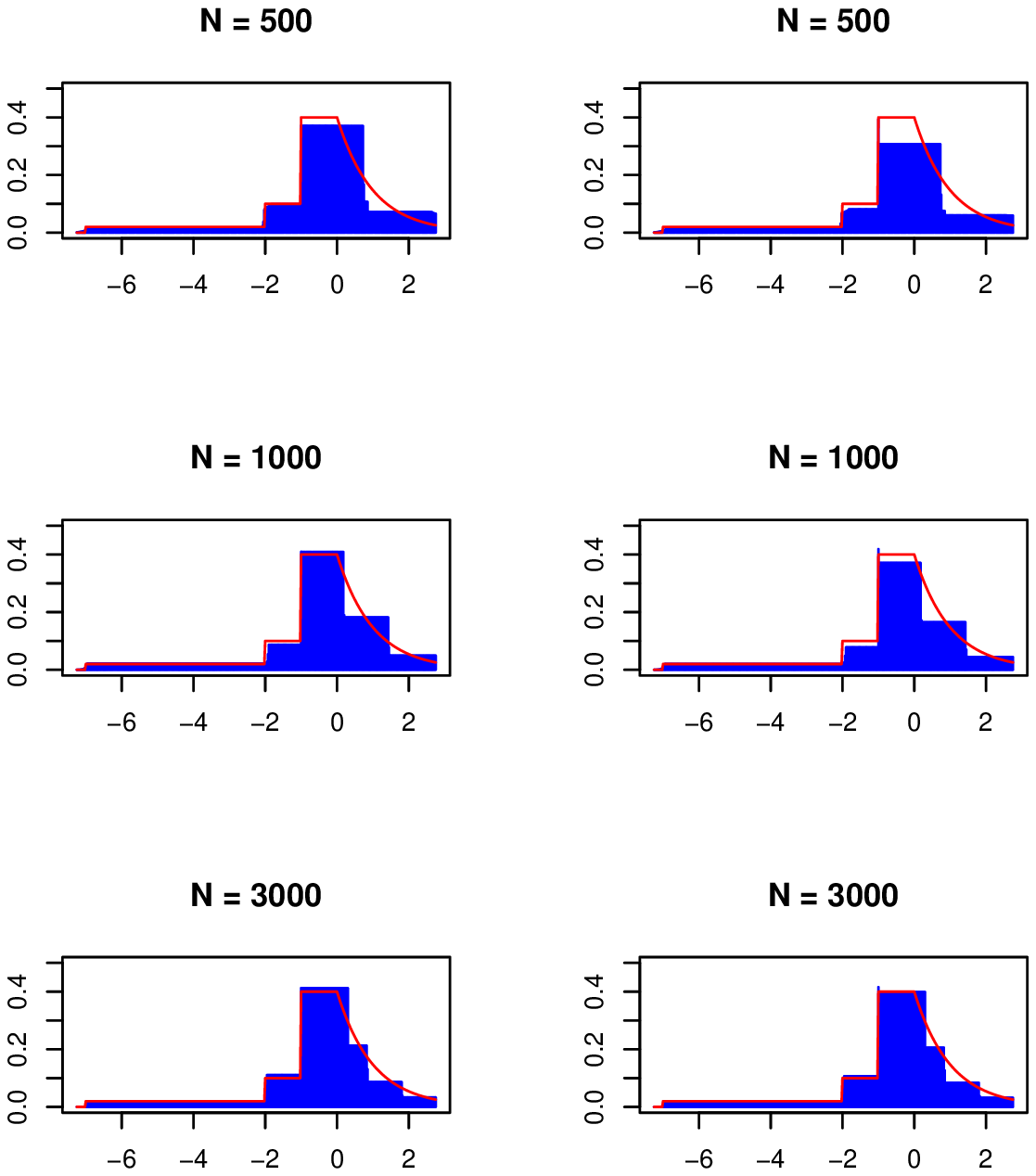}\\
\renewcommand{\baselinestretch}{1}
  \parbox{5in}
  {\caption{\quad The true unimodal density $\lambda_1(t)$~(solid line) and
  weighted average density estimates given the mode $\theta = -1$
  produced by the SIP sampler~(Steps~(ii) and~(iii)
  of Algorithm~\ref{bigSIS}) based on $a=0$ and $b=1$~(left column)
  and $a=0.9$ and $b=100$~(right column) for $G$.}
  \label{fig:SIP1}}
\end{figure}
\begin{figure}[t]
\centering\includegraphics[width=6in]{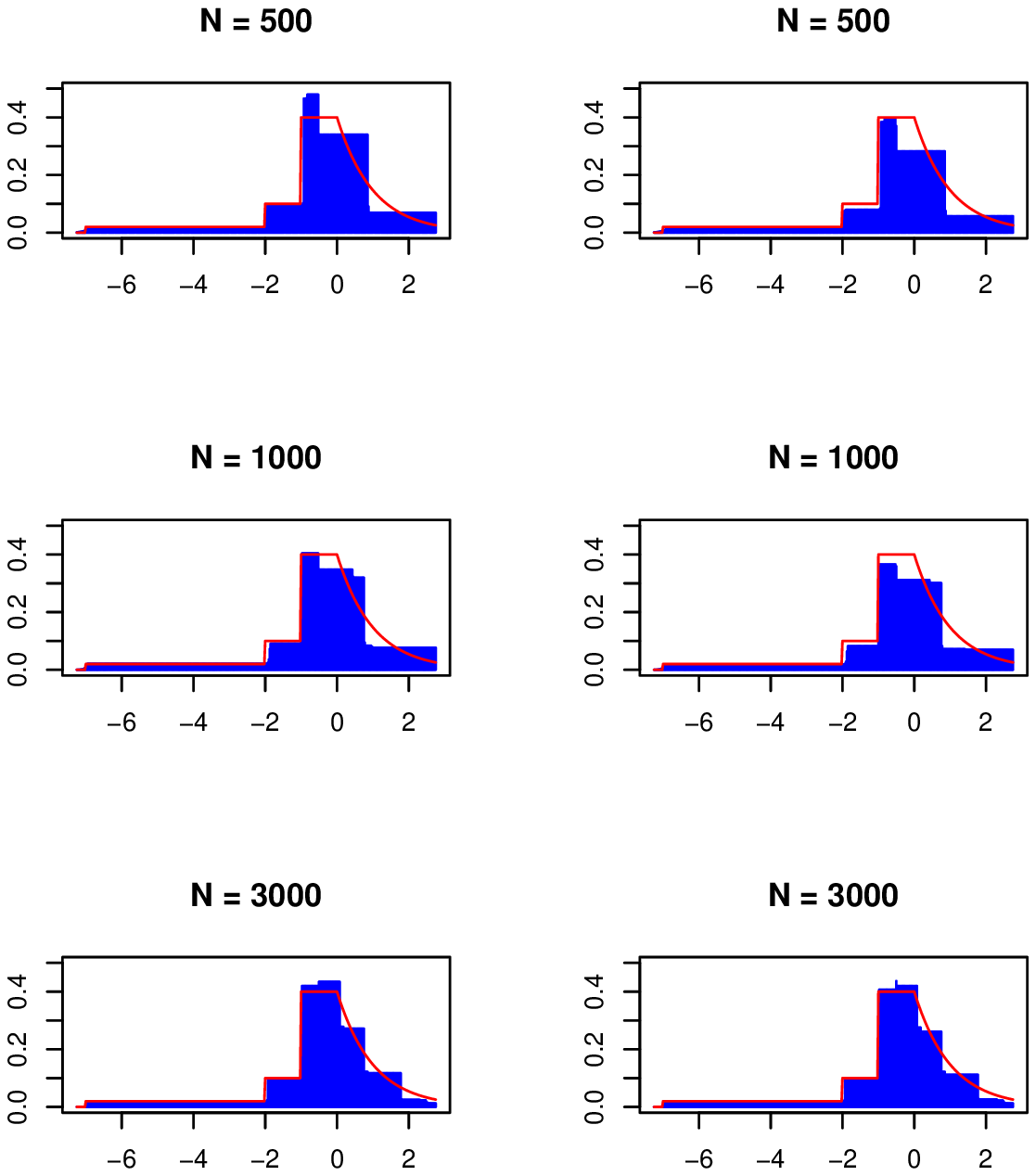}\\
\renewcommand{\baselinestretch}{1}
  \parbox{5in}
  {\caption{\quad The true unimodal density $\lambda_1(t)$~(solid line) and
  weighted average density estimates given the mode $\theta = -0.5$
  produced by the SIP sampler~(Steps~(ii) and~(iii)
  of Algorithm~\ref{bigSIS}) based on $a=0$ and $b=1$~(left column)
  and $a=0.9$ and $b=100$~(right column) for $G$.}
  \label{fig:SIP2}}
\end{figure}
\begin{figure}[t]
\centering\includegraphics[width=6in]{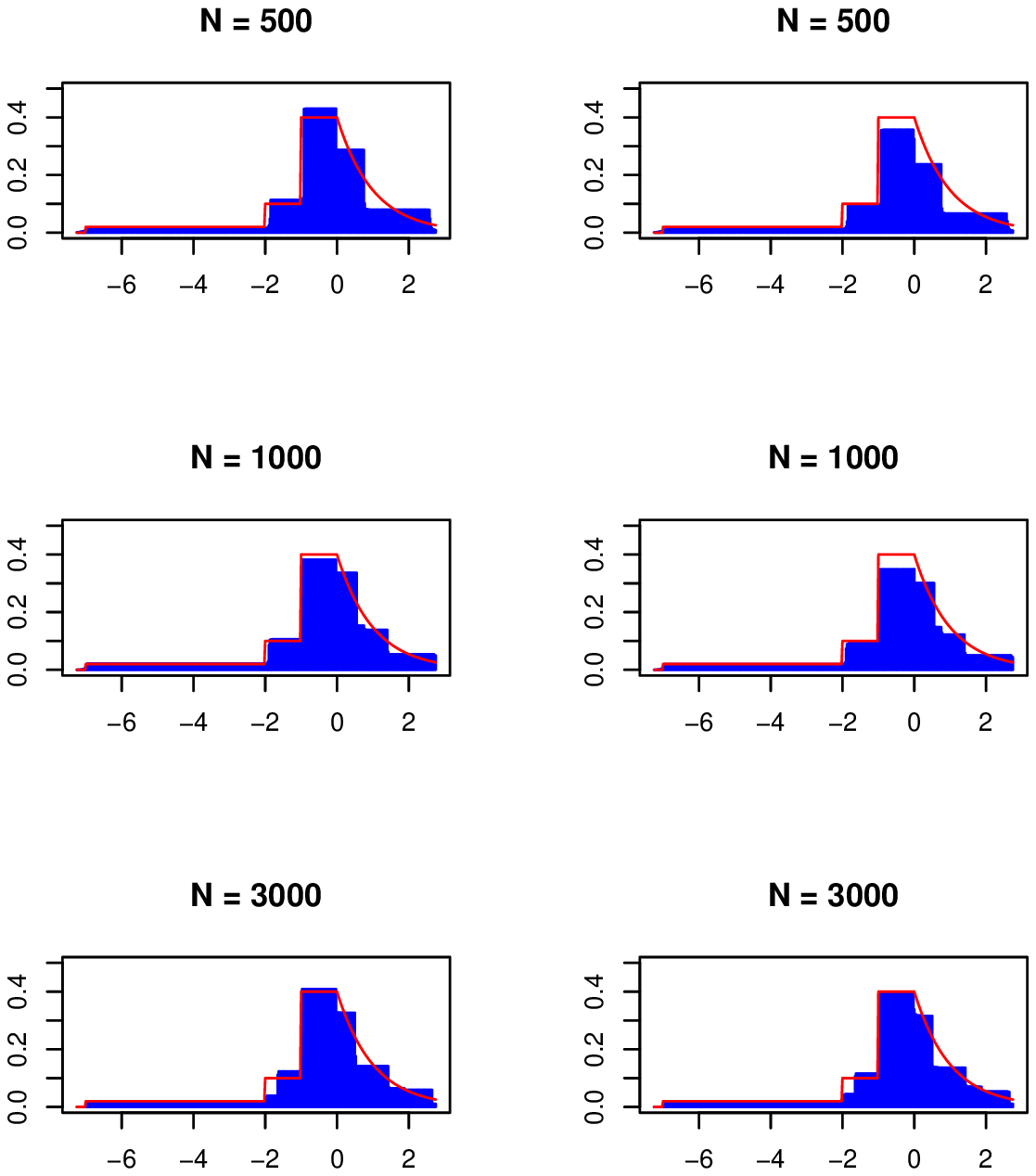}\\
\renewcommand{\baselinestretch}{1}
  \parbox{5in}
  {\caption{\quad The true unimodal density $\lambda_1(t)$~(solid line) and
  weighted average density estimates given the mode $\theta = 0$
  produced by the SIP sampler~(Steps~(ii) and~(iii)
  of Algorithm~\ref{bigSIS}) based on $a=0$ and $b=1$~(left column)
  and $a=0.9$ and $b=100$~(right column) for $G$.}
  \label{fig:SIP3}}
\end{figure}
\begin{figure}[t]
\centering\includegraphics[width=6in]{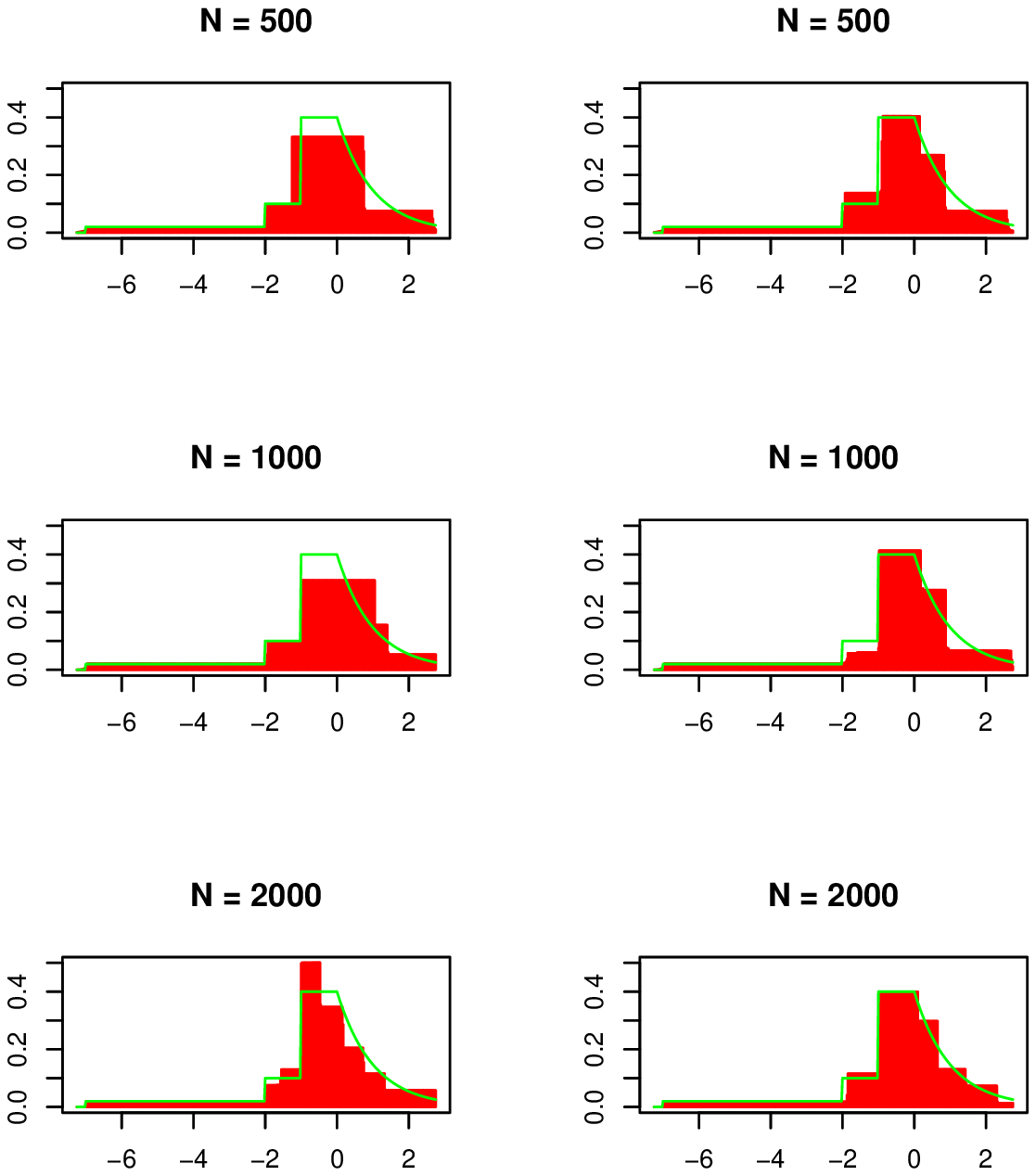}\\
\renewcommand{\baselinestretch}{1}
  \parbox{5in}
  {\caption{\quad The true unimodal densities $\lambda_1(t)$~(solid lines) and
  weighted average density estimates
  produced by Algorithm~\ref{bigSIS} based on a $N(0,1)$
  density~(left column) and a
  $N(0,0.25^2)$ density~(right column) for $\rho(\theta)$.}
  \label{fig:bigSIS1}}
\end{figure}
\begin{figure}[t]
\centering\includegraphics[width=6in]{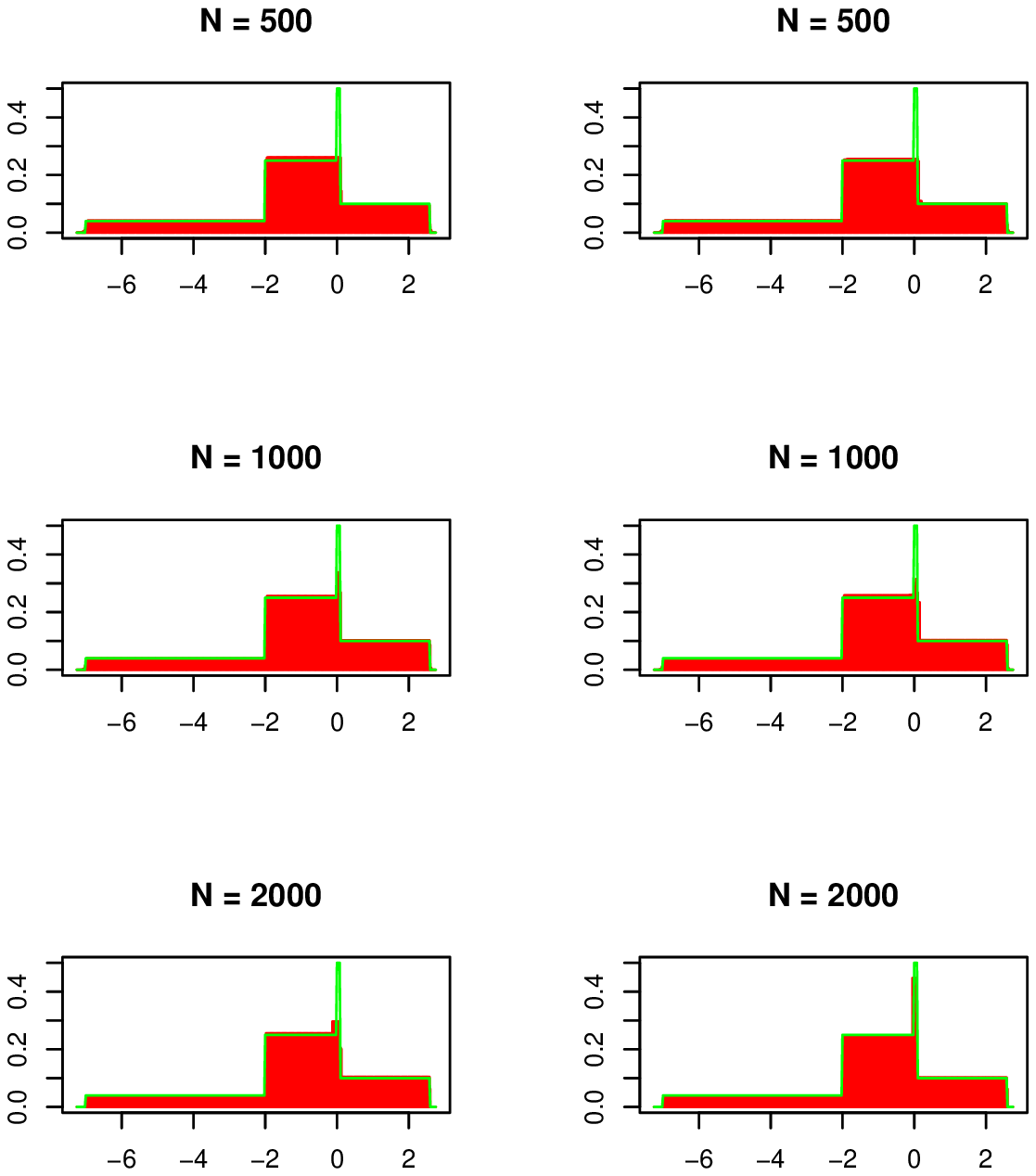}\\
\renewcommand{\baselinestretch}{1}
  \parbox{5in}
  {\caption{\quad The true unimodal densities $\lambda_2(t)$~(solid lines) and
  weighted average density estimates
  produced by Algorithm~\ref{bigSIS} based on a $N(0,1)$
  density~(left column) and a
  $N(0,0.25^2)$ density~(right column) for $\rho(\theta)$.}
  \label{fig:bigSIS2}}
\end{figure}
\begin{figure}[t]
\centering\includegraphics[width=6in]{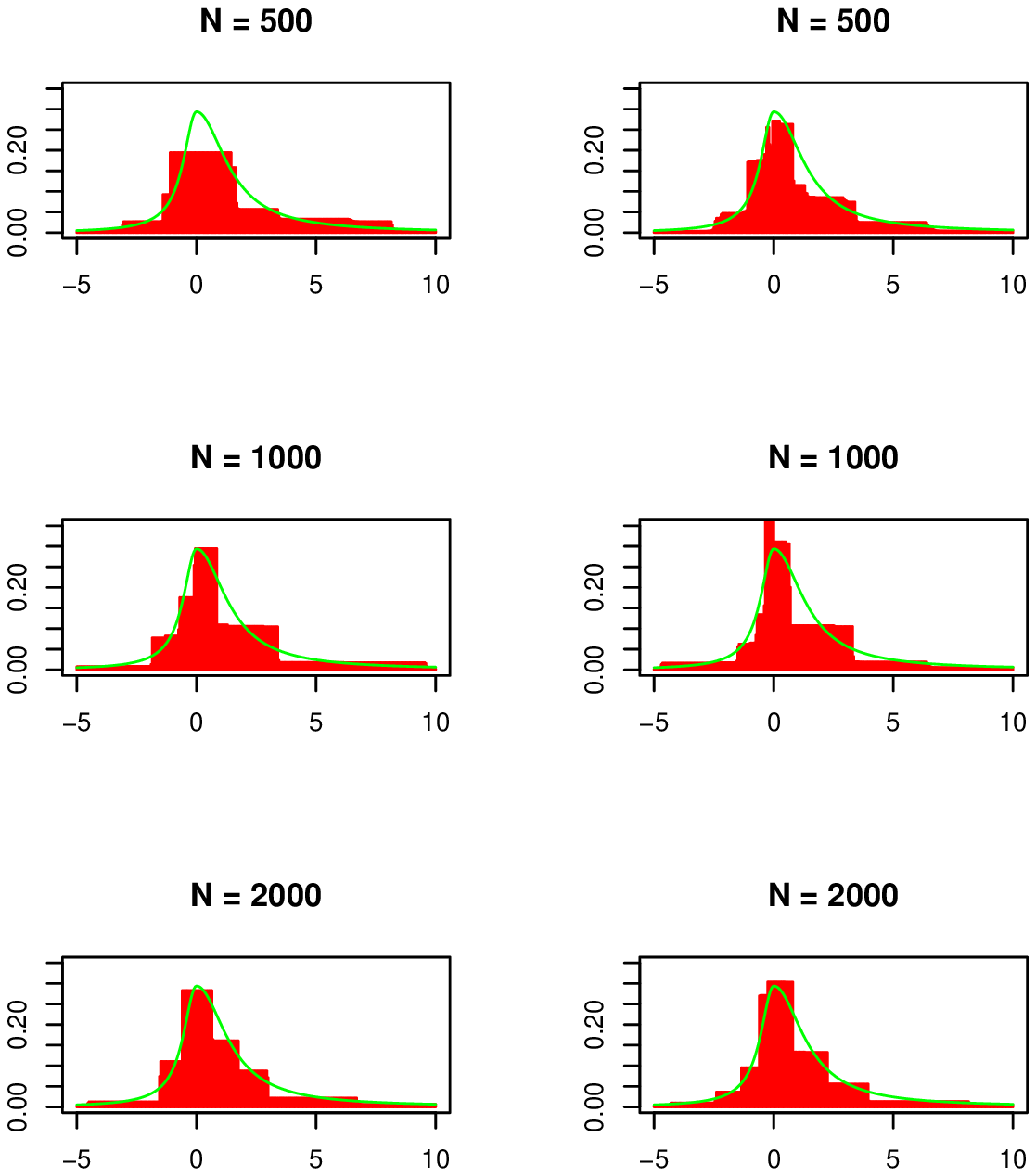}\\
\renewcommand{\baselinestretch}{1}
  \parbox{5in}
  {\caption{\quad The true unimodal densities $\lambda_3(t)$~(solid lines) and
  weighted average density estimates
  produced by Algorithm~\ref{bigSIS} based on a $N(0,1)$
  density~(left column) and a
  $N(0,0.25^2)$ density~(right column) for $\rho(\theta)$.}
  \label{fig:bigSIS3}}
\end{figure}

\end{document}